\documentclass[preprint,3p]{elsarticle} 

\usepackage{amssymb}
\usepackage{lineno,hyperref,amsmath,mathtools}
\usepackage{bm}
\usepackage{xcolor}
\usepackage{subcaption}

\usepackage[utf8]{inputenc}
\usepackage[noline,boxed]{algorithm2e}
\usepackage[skins]{tcolorbox}
\usepackage{algorithmic}
\usepackage{mathrsfs}
\usepackage{float}
\usepackage{tabularx} 

\usepackage{soul}

\newcommand{\ve}[1]{\textbf{\textit{#1}}}

\newcolumntype{b}{X}
\newcolumntype{s}{p{0.2cm}}
\newcolumntype{t}{p{0.45cm}}

\begin{document}

\begin{frontmatter}



\title{An improved local radial basis function method for solving small-strain elasto-plasticity}


\affiliation[fs]{organization={Faculty of Mechanical Engineering, University of Ljubljana},
	addressline={Aškerčeva cesta 6}, 
	city={Ljubljana},
	postcode={SI-1000}, 
	country={Slovenia}}
\affiliation[imt]{organization={Institute of Metals and Technology},
	addressline={Lepi pot 11}, 
	city={Ljubljana},
	postcode={SI-1000}, 
	country={Slovenia}}

\author[fs]{Gašper Vuga}
\ead{gasper.vuga@fs.uni-lj.si}
\author[fs,imt]{Boštjan Mavrič}
\ead{bostjan.mavric@it.uu.se}
\author[fs,imt]{Božidar Šarler\corref{corr}}
\ead{bozidar.sarler@fs.uni-lj.si}
\cortext[corr]{Corresponding author}    

%
%

\begin{abstract}
Strong-form meshless methods received much attention in recent years and are being extensively researched and applied to a wide range of problems in science and engineering. However, the solution of elasto-plastic problems has proven to be elusive because of often non-smooth constitutive relations between stress and strain. The novelty in tackling them is the introduction of virtual finite difference stencils to formulate a hybrid radial basis function generated finite difference (RBF-FD) method, which is used to solve small-strain von Mises elasto-plasticity for the first time by this original approach. The paper further contrasts the new method to two alternative legacy RBF-FD approaches, which fail when applied to this class of problems. The three approaches differ in the discretization of the divergence operator found in the balance equation that acts on the non-smooth stress field. Additionally, an innovative stabilization technique is employed to stabilize boundary conditions and is shown to be essential for any of the approaches to converge successfully. Approaches are assessed on elastic and elasto-plastic benchmarks where admissible ranges of newly introduced free parameters are studied regarding stability, accuracy, and convergence rate.

\end{abstract}  




\end{frontmatter}


\section{Introduction} \label{sec: Intro}
Computational modeling of elasto-plastic solid mechanics is a well-established field of science \cite{argyris_elasto-plastic_1965}. Most of the research has been, since the mid-60s, performed using mesh-based numerical methods to solve the related governing partial differential equations (PDEs). These methods, in general, fall into three main groups, i.e., finite element methods (FEM) \cite{de_souza_neto_computational_2011}, finite volume methods (FVM) \cite{cardiff_thirty_2021}, and boundary element methods (BEM) \cite{katsikadelis_boundary_2016}. The most widely used method for solving continuum mechanics problems is FEM, which many commercial packages adapt. The main drawback of mesh-based methods is the time-consuming quality mesh generation - polygonization.

As an alternative to mesh-based methods, meshless methods (MMs) \cite{atluri_meshless_2004, liu_introduction_2005, liu_meshfree_2009, sarler_recent_2010, li_meshless_2013, pepper_introduction_2014} began to evolve. In MMs, the geometry is discretized by a cloud of nodes distributed over the domain of interest. This allows flexibility when dealing with complex geometries, moving boundary problems and problems with multiple dimensions. Also, various discretization adaptivity types can be achieved \cite{slak_adaptive_2020}, and multilevel techniques can easily be applied \cite{zamolo_novel_2019}. 

Many different MMs have been developed so far, such as the element-free Galerkin method (EFGM) \cite{belytschko_element-free_1995}, local Petrov Galerkin method (LPGM) \cite{atluri_new_1998}, direct meshless local Petrov Galerkin method (DMLPGM) \cite{mirzaei_direct_2016}, smooth particle hydrodynamics method (SPH) \cite{gingold_smoothed_1977}, reproducing kernel particle method (RKPM) \cite{liu_reproducing_1995}, finite point method (FPM) \cite{onate_finite_1996}, radial point interpolation method (RPIM) \cite{wang_point_2002}, method of fundamental solutions (MFS)  \cite{liu_method_2019}, diffuse approximate method (DAM) \cite{prax_collocated_1996}, local radial basis function collocation method (LRBFCM) \cite{sarler_meshfree_2006}, known in recent years also as radial basis generated finite differences (RBF-FD) \cite{tolstykh_using_2003} etc.

MMs can be generally classified according to their formulation procedure (weak, strong, or combined formulation), function approximation (moving and weighted least-squares approximation, integral representation method or point interpolation method) and domain representation (domain-type or boundary-type discretization method) \cite{liu_introduction_2005}. The studies with EFGM \cite{kargarnovin_elasto-plastic_2004}, LPGM \cite{atluri_new_1998}  and RKPM \cite{chen_reproducing_1996, ji-fa_meshfree_2005} on solving elasto-plastic and other non-linear mechanical problems have successfully demonstrated the use of domain-type weak-form MMs.

Weak-form MMs possess good stability and accuracy and naturally satisfy the Neumann boundary conditions (BCs). In contrast to FEM, a high-order convergence and high resolution of steep gradients can be efficiently achieved \cite{lu_new_1994}. However, numerical integration makes weak-form MMs computationally expensive. Also, the background mesh (or so-called shadow mesh \cite{atluri_new_1998}), used to evaluate domain integrals, must be constructed, which again demands some kind of meshing \cite{liu_meshfree_2003}. As opposed to the weak-form MMs, in strong-form MMs, PDEs are directly discretized without the need for a weak (integral) form. Implementation is also easier, and they are computationally efficient and do not include any background mesh for neither approximation nor integration. Nonetheless, proving the stability of such methods still remains an open problem, especially when Neumann BCs are present \cite{lu_new_1994}.

This paper addresses the possibility of solving two-dimensional elasto-plastic solid mechanics with strong-form MM.  The RBF-FD method \cite{zhang_meshless_2000, tolstykh_using_2003} (also LRBFCM, \cite{sarler_meshfree_2006}) is employed, and the approximation functions, radial basis functions (RBFs), are chosen to be polyharmonic splines (PHSs). RBF-FD method is a local collocation method where PDEs are discretized in a strong form. It can be interpreted as a generalization of the finite difference (FD) method \cite{ozisik_finite_2017} and has its roots in the original Kansa method, where collocation is performed globally \cite{kansa_multiquadricsscattered_1990-1, kansa_multiquadricsscattered_1990}. In the Kansa method, the resultant stiffness matrix is large, dense and ill-conditioned, making it a computationally very challenging method. In RBF-FD, the solution is locally approximated on a local support domains. Local approximation of the solution field results in small full systems that have to be inverted for each node in order to obtain local weights of differential operators (DOs). With all weights computed, the global stiffness matrix can be composed. Because of its sparse structure, it is better conditioned and easy to solve than the global Kansa method.


RBF-FD/LRBFCM method has been successfully applied to various academic and industrial problems such as diffusion \cite{sarler_meshfree_2006}, diffusion-convection with phase change \cite{vertnik_meshless_2006}, natural convection \cite{kosec_solution_2013}, macrosegregation \cite{kosec_simulation_2014}, radiation-diffusion on surfaces \cite{lehto_radial_2017}, natural convection under the influence of a static magnetic field \cite{mramor_simulation_2013, mramor_application_2020}, flow in porous media \cite{hatic_meshless_2020}, the flow of a non-Newtonian fluid \cite{hatic_meshless_2021},  phase field modeling of dendritic solidification \cite{dobravec_reduction_2020}, micro-combustion \cite{bayona_micro-combustion_2021}, transient direct-chill aluminum billet casting problem with simultaneous material and interphase moving boundaries \cite{vertnik_solution_2006} to name a few.

The first use of RBF-FD in solid mechanics was presented in \cite{zhang_meshless_2000} on two standard 2D linear-elastic benchmarks, where different types of RBFs were tested. In \cite{tolstykh_using_2003}, the first use of Hardy multiquadrics (MQs) as RBFs was presented on various linear elastic cases. The MQs approach is still being used nowadays \cite{ferreira_analysis_2003, ferreira_computation_2006, stevens_solution_2013, hanoglu_hot_2019}. Solutions for thermo-elasticity problems were initially presented in \cite{gerace_localized_2006} and extended in \cite{mavric_collocation_2014} for coupled problems.

The first attempt of solving a non-linear mechanical problem was presented in \cite{hanoglu_thermo-mechanical_2011} where a hot-rolling process of steel bars was modeled. The updated Lagrangian formulation in 2D has been employed where MQs were used for RBFs. The material was treated as ideally plastic, and a non-linear system of equations was solved via direct iteration. The BCs were included in local interpolation problems. The proposed approach has also been successfully used in a multi-pass hot-rolling with non-symmetric groove types \cite{hanoglu_multi-pass_2018, hanoglu_rolling_2018, hanoglu_hot_2019}. In \cite{mavric_meshless_2017}, a more complicated elasto-visco-plastic constitutive model was used for modeling the stress-strain state in solidified part the during direct-chill casting of aluminum billets. Simulations have been performed in axisymmetric coordinates, and once again, the MQs were employed as RBFs. During the iterative process of solving a system of non-linear equations, the Jacobian was determined numerically, resulting in slow convergence. With the same solution procedure, a 2D slice model of continuous casting of steel billets \cite{mavric_meshless_2020} has been successfully implemented. In \cite{strnisa_meshless_2022}, the first use of RBF-FD with PHS was presented for solving elasto-plastic problems. Only elastic material parameters were used when composing Jacobian, so many global iterations were needed for the convergence.    


This work presents the derivation and implementation of a general elasto-plastic solver, where the RBF-FD method is used for discretizing a governing PDE and a return-mapping algorithm (RMA) is used for solving constitutive relations and determining consistent tangent operator (CTO) at the material point level. Von Mises yield function is used, and hardening is set to be isotropic. The system of non-linear equations is solved by the Newton-Raphson iteration algorithm (NRIA). PHSs \cite{flyer_role_2016, bayona_role_2017} are used as the RBFs since they do not include an undetermined shape parameter that must be set when using other RBFs such as multiquadrics, inverse multiquadrics, Gaussians and others. The process of obtaining optimal shape parameter can be computationally very intensive, which is not desired when dealing with large engineering problems. The use of PHSs has gained popularity in recent years and has been successfully applied to non-linear elliptic problems \cite{bayona_role_2017}, solution of dendritic solidification \cite{dobravec_reduction_2020} and linear elastic problems \cite{slak_adaptive_2020}.      

This work introduces and assesses three different meshless discretization approaches of the non-linear boundary value problem (BVP). These approaches are here denoted \textit{direct}, \textit{composed} and \textit{hybrid}. With a \textit{direct} approach, the divergence operator in the balance equation is applied to the stress field in its continuous form, so the derivatives up to the second order are discretized when composing a Jacobian. This kind of RBF-FD discretization has been successfully used on linear problems \cite{slak_adaptive_2020, mavric_collocation_2014} and previously mentioned hot-rolling of steel bars (ideal-plasticity) \cite{hanoglu_thermo-mechanical_2011, hanoglu_multi-pass_2018, hanoglu_rolling_2018, hanoglu_hot_2019}, and direct-chill of aluminum billets (elasto-visco-plasticity) \cite{mavric_meshless_2017}. In the \textit{composed} approach, first presented in the present paper, the divergence operator is discretized and acts on the discretized form of the stress field. Derivatives only up to the first order are computed. In the \textit{hybrid} approach, secondary variables (stresses and strains) are discretized on a new set of virtual nodes - secondary nodes (SNs). To each collocation node (CN) a different FD stencil is prescribed where points of FD stencils coincide with SNs. The divergence operator is then discretized via the FD method. A similar approach was used for solving 3D thermal problems \cite{gerace_model-integrated_2014} and compressible fluid flow \cite{harris_application_2015}. In studies, the field variables were RBF-interpolated on a $2^{\text{nd}}$ order virtual FD stencils adapted to each CN separately. Then all the necessary differential operators were determined via the FD method. Neither of these studies used PHSs. As mentioned initially, the strong-from MMs can have problems with Neumann BCs. This paper proposes a new simple but crude stabilization technique, which turns out to be essential for the convergence of the method.

The main originality of the present paper is an introduction and a verification of \textit{direct}, \textit{composed} and \textit{hybrid} approaches used for the discretization of governing equations of elasto-plasticity, with the use of RMA.
It is shown that the \textit{hybrid} approach significantly outperforms the other two approaches in terms of stability, accuracy and convergence. This work represents the first attempt at applying a strong form MM for solving elasto-plasticity problems, where Jacobian is determined using a consistent tangent operator.


The present paper is structured in the following way; first, the governing equations of solid mechanics are introduced in Section \ref{sec: GOVeq}, where special attention is given to a constitutive model associated with von Mises elasto-plasticity. In Section \ref{sec: NUM_MET}, a presentation of the RBF-FD method is given. Firstly, the geometry discretization and local support construction for local approximation is described. Then, the procedure of constructing a local interpolant and its application on the approximation of differential operators with the prescription of radial basis function is presented. Next, the procedure of solving a non-linear system of equations via the NRIA with originally introduced spatial discretization is described. A presentation of a new approach of BC stabilization and a brief revision of the integration algorithm of elasto-plastic constitutive relations is given at the end of this section. In Section \ref{sec: NUM_EX_EL}, a parametric study on elastic benchmarks is performed where the effect of newly introduced parameters is assessed. The findings are then used for calculations of elasto-plastic case of the internally pressurized annulus, and a comparison between the three approaches is given in \ref{sec: NUM_EX_PL}. The conclusions are presented in Section \ref{sec: Conc}.

\section{Governing equations}\label{sec: GOVeq}
In this work, the equations of small strain elasto-plasticity with associated von Mises flow rule and isotropic hardening are being solved. Although plane strain and plane stress cases are studied in this work, the equations are presented in a coordinate-free form.

Consider the continuum material occupying the domain $\Omega$, with boundary $\Gamma$. For each point within $\Omega$, the equilibrium state is described with a balance law
\begin{equation}\label{eq: balance DE}
	\nabla \cdot \bm\sigma  = - \ve{f}, 
\end{equation}
where $\bm\sigma$ is the Cauchy stress tensor and $\ve{f}$ is the body force vector. Strain tensor is in the case of small deformations defined as
\begin{equation}\label{eq: TSDeformacij}
	\bm\varepsilon = \nabla^s \ve{u}, \qquad  \nabla^s = (\nabla + \nabla^\top)/2 ,
\end{equation}
where $\nabla^s$ denotes the symmetric gradient operator and $\ve{u}$ the displacement vector. In order to obtain a unique solution of the equilibrium equation (\ref{eq: balance DE}) in terms of displacement vector $\ve{u}$, appropriate BCs must be prescribed at the boundary $\Gamma = \Gamma_u \cup  \Gamma_T \cup \Gamma_F$ :
\begin{equation}\label{eq: robni pogoji}
	\begin{aligned}
		&   \ve{u} = \bar{\ve{u}}     &&  \text{on} \quad \Gamma_u, \\
		&   \bm\sigma \cdot \ve{n} = \bar{\ve{T}}   &&  \text{on} \quad \Gamma_T, \\
		&   \left\{u_n, T_t\right\}  = \left\{0,0\right\}                  &&  \text{on} \quad \Gamma_F, \\
	\end{aligned}
\end{equation}
where $\bar{\ve{u}}$ is the prescribed value of $\ve{u}$ on $\Gamma_u$. On $\Gamma_T$ the traction vector $\bar{\ve{T}}$ is prescribed, which is equal to the stress projected in a normal direction $\ve{n}$. On free surfaces, its value is zero. On free-slip boundary $\Gamma_F$, displacement in normal direction $u_n = \ve{u} \cdot \ve{n}$ and traction in tangential direction $T_t = \ve{T} \cdot \ve{t}$ are prescribed, both equal to zero. The material is at $\Gamma_F$ restricted to move in the normal direction and free in the tangential direction, with no friction imposed. The kind of BC where $\Gamma_F$ is a straight line is called symmetry BC.

The relationship between stress and strain is defined by a constitutive law.
In case the material is fully recoverable, stress can be directly determined from strain using Hooke's constitutive law
\begin{equation} \label{eq: Hooks Law}
	\bm\sigma = \textbf{\textsf{D}}^e : \bm\varepsilon^e,
\end{equation}
where $\textbf{\textsf{D}}^e$ represent the fourth-order elasticity tensor and $\bm\varepsilon^e$ elastic strain tensor. For isotropic elastic material, $\textbf{\textsf{D}}^e$ can be defined by only two material parameters, and equation (\ref{eq: Hooks Law}) can be written as
\begin{equation}
	\bm\sigma = 2 \mu \bm\varepsilon^e + \lambda \bm{I} \text{tr}(\bm\varepsilon^e), 
\end{equation}
where $\bm{I}$ is the second-order identity tensor. Parameters $\mu$ and $\lambda$ are Lam\'e constants that can be expressed with experimentally measured Young's modulus $E$ and Poisson's ratio $\nu$ as $\mu = E/(2(1+\nu))$ and $\lambda = \nu E /((1+\nu)(1-2\nu))$. The parameter $\mu$ is also known as the shear modulus $G$.

Once the material reaches the elastic limit or the so-called yield stress, it begins to yield, and the plastic strain $\bm\varepsilon^p$ is induced. By assuming small deformations, the total strain can be additively split as
\begin{equation}\label{eq: additive_split}
	\bm\varepsilon = \bm\varepsilon^e+\bm\varepsilon^p.
\end{equation}
Using this relation, Hooke's law can be rewritten as
\begin{equation}
	\bm\sigma = \textbf{\textsf{D}}^e : (\bm\varepsilon - \bm\varepsilon^p).
\end{equation}
The yield criterion defines a critical stress state at which material starts yielding. Using the von Mises yield criterion, the plastic flow will occur when the von Mises equivalent stress $\sigma_{vm} = \sqrt{3 J_2(\bm{s})}$ equals the uniaxial yield stress $\sigma_y$, obtained from the uniaxial tensile test. Here $J_2$ represents the second invariant of a deviatoric stress tensor  $\ve{s} = \bm\sigma - \ve{I}\text{tr}(\bm\sigma)/3$. Uniaxial yield stress is generally a non-linear function of the hardening thermodynamic force $\sigma_y = \sigma_y(\tau)$. Assuming isotropic hardening, $\tau$ becomes a function of a single scalar value $\bar{\varepsilon}^p$ called the accumulated plastic strain. The hardening curve is typically specified in the form $\sigma_y(\kappa) = \sigma_{y0}  + \tau(\bar{\varepsilon}^p)$, where $\sigma_{y0}$ stands for the initial yield stress of virgin material.

The admissible region in stress space can be compactly specified with the yield function $\Phi(\bm\sigma, \tau)$ where $\Phi = 0$ when the material is subjected to plastic flow and $\Phi < 0$ when the response is elastic. Using the von Mises yield criterion; it is defined as
\begin{equation}\label{eq: vonMises YF}
	\Phi(\bm\sigma, \tau) = \sqrt{3 J_2(\bm{s}(\bm\sigma))} - \sigma_y(\tau), \qquad \tau = \tau(\bar{\varepsilon}^p).
\end{equation}
Next, the plastic flow rule defines how the plastic strain evolves. Within the plastic potential flow theory, it is assumed that $\bm\varepsilon^p$ is evolving in the direction of a plastic flow vector $\bm{N} = \partial \Psi / \partial \bm\sigma$, where  $\Psi$ represents plastic flow potential. For ductile materials, $\Psi = \Phi$ is assumed. Then, the so-called associated plastic flow rule is, in the rate form, defined as
\begin{equation}
	\dot{\bm\varepsilon}^p = \dot{ \gamma} \bm{N} = \dot{ \gamma} \frac{\partial \Phi}{\partial \bm\sigma} = \dot{ \gamma}\sqrt{\frac{3}{2}} \frac{\bm{s}}{||\bm{s}||},
\end{equation}
where the von Mises yield function is substituted from equation (\ref{eq: vonMises YF}) and $||\bm{s}|| = \sqrt{\bm{s} : \bm{s}}$. In this case, the plastic flow vector $\bm{N}$ is called the Prandtl-Reuss flow vector. A non-negative plastic multiplier $\gamma$ represents the magnitude of the plastic strain and connects the stress space with the plastic strain space.

The evolution of the hardening internal state variable, which is, in this case, equivalent to $\bar{\varepsilon}^p$, is within the previous assumptions given as 
\begin{equation}\label{eq: Hardening Law}
	\dot{\bar{\varepsilon}}^p = \dot{ \gamma}, \qquad  \quad \bar{\varepsilon}^p = \sqrt{\frac{2}{3}} ||\bm\varepsilon^p||.
\end{equation}
Lastly, a set of loading/unloading (also known as Kuhn-Tucker) conditions that specify when the evolution of plastic strain and internal variables may occur is given as
\begin{equation}\label{eq: K-T}
	\Phi(\bm\sigma, \tau) \leq 0, \quad \dot{\gamma} \geq 0, \quad \Phi(\bm\sigma, \tau) \dot{\gamma}=0.
\end{equation}
In-detail description of the above relations can be found in \cite{de_souza_neto_computational_2011, simo_computational_2000}.

\section{Numerical method} \label{sec: NUM_MET}

\subsection{RBF-FD}
RBF-FD method can be represented as a generalization of the classical finite difference (FD) method, where a spatial differential operator (DO) at a node is approximated as the weighted sum of the field values. Unlike the FD method, where DOs are computed on regular node arrangements (RNAs), the RBF-FD method allows an evaluation of DOs also on scattered node arrangements (SNAs). This makes problems on complex geometries solvable but comes at the price of calculating weights for each discretization node separately. Compared with the FD method, one of the main advantages of RBF-FD is simple local refinement in critical areas. The number of weights that belong to a specific node depends on the number of neighboring support nodes enclosed inside a local support domain. A higher number of these nodes means better approximation but also less efficient calculation of the local weights. Also, the global sparse stiffness matrix becomes denser when using large supports, which can cause problems when solving a system of equations and saving non-zero values for large cases. 


\subsubsection{Geometry discretization in 2D}\label{sec: Geom_dist}
\begin{figure}
	\begin{centering}
		\includegraphics[width=\linewidth]{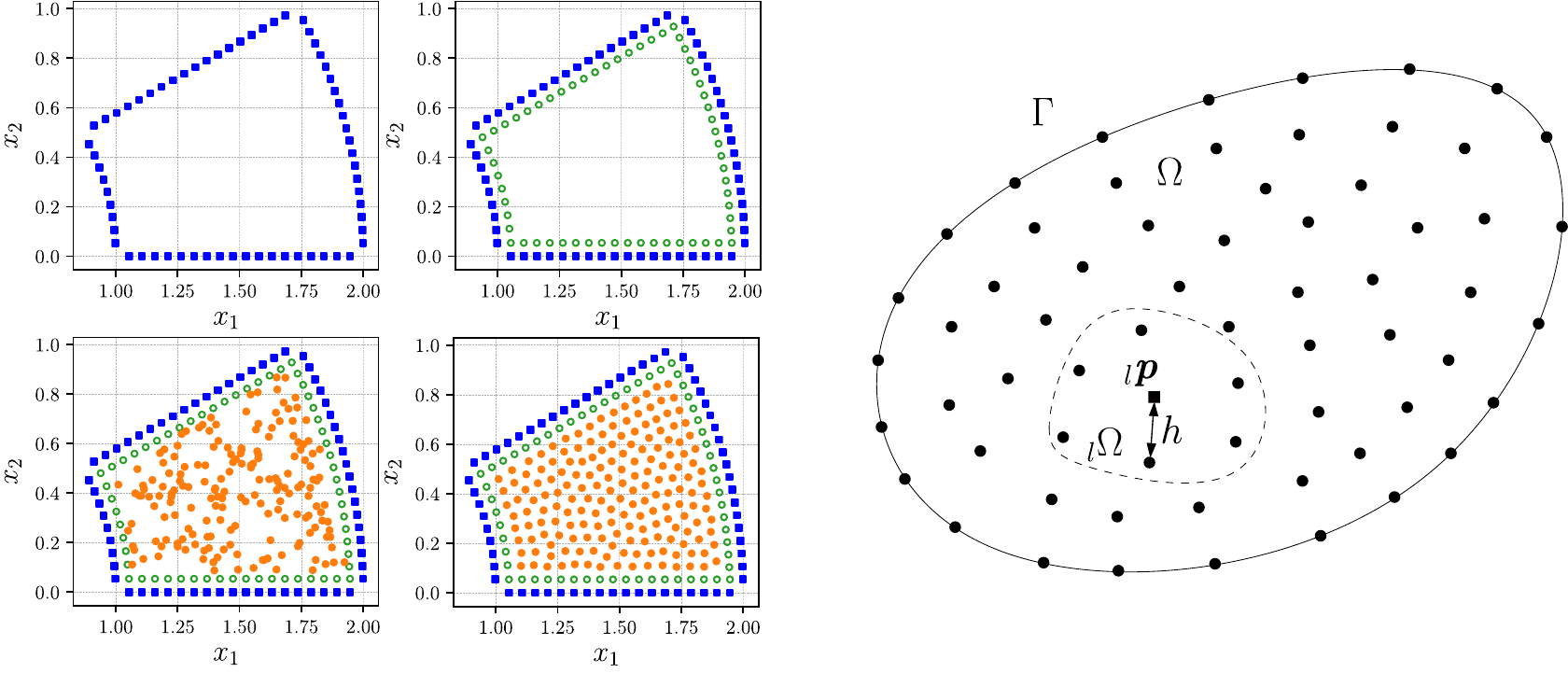}
		\caption{Left: four steps of nodes positioning on the annulus section geometry. Boundary nodes are represented with squares (top-left), the addition of the inner boundary nodes with hollow circles (top-right), the addition of randomly distributed inner nodes with full circles (bottom-left) and homogeneous redistribution of inner nodes (bottom-right). Right: discretization of generic domain $\Omega$ with its boundary $\Gamma$. The square represents the central node ${}_{l}\ve{p}$ of a local support domain ${}_{l}\Omega$ that is defined with ${}_{l}N$ support nodes inside a dashed boundary. $h$ measures the distance to the closest neighbor.}
		\label{fig: combined}
	\end{centering}
\end{figure}

The presented algorithm defines plane geometry as a set of parameterized curves with corresponding functions of outward-facing normal vectors. Discretization with SNA starts with prescribing node density $\rho(\ve{p})$ over the investigated domain $\Omega$ including the boundary $\Gamma$, where $\ve{p}$ represents the position vector. From here, the distance between two consecutive nodes is calculated as $h(\ve{p}) = \theta/\sqrt{\rho(\ve{p})}$, where $\theta$ represents geometric factor defined for an ideal hexagonal lattice as $\theta = \sqrt{2/\sqrt{3}}$. Knowing $h(\ve{p})$ and the prescribed geometry, $N_{b}$ number of boundary nodes are first positioned in a manner that can be observed in Figure \ref{fig: combined} (left), on the example of an annulus section. The corner nodes are omitted since if different types of BCs are prescribed on the sides that meet at the corner, the question of which BC type should be prescribed to a corner node arises. In terms of constructing a local interpolant (described in the following Section \ref{sec: Interpolant_construction}), where the evaluation point lies on the boundary of the local support it can result in bad conditioning of the interpolant. After positioning boundary nodes, inner boundary nodes are positioned in the opposite direction of the boundary normal vectors (hollow circles). Adding this inner boundary improves local interpolation at boundaries where Neumann BCs ($\Gamma_T$ and $\Gamma_F$) are prescribed \cite{mavric_collocation_2014}. Next, inner nodes (full circles) are randomly added inside a domain defined by the inner boundary. Those are then homogeneously redistributed with a minimization process that mimics electrostatic repulsion. This is done to achieve local isotropy of node distribution. The elaboration of node generation can be found in \cite{mavric_meshless_2017}.

\subsubsection{Local support domain construction}\label{sec: Sub_cons}
In order to locally approximate an arbitrary function, a local support domain ${}_{l}\Omega$, where $l$ runs over all discretization nodes $l = 1, ..., N$, should be defined. An illustration of some specific ${}_{l}\Omega$, centered at node ${}_{l}\ve{p}$ is shown in Figure \ref{fig: combined} (right). The choice of ${}_{l}N$ number of support nodes that defines ${}_{l}\Omega$, is based on prescribing the ${}_{l}N - 1$ nearest nodes to the central node. If the node arrangement is homogeneously isotropic, the choice of nearest neighbors provides good support domains for approximation. This holds for nodes not on or near the boundary, where some other conditions of choosing ${}_{l}N$ have to be prescribed (elaboration can be found in \cite{davydov_adaptive_2011, mavric_meshless_2017}). When ${}_{l}\ve{p} \in \Gamma$, there is no restriction imposed on the maximum number of boundary nodes inside ${}_{l}\Omega$, in contrast with previous works \cite{hanoglu_thermo-mechanical_2011, hanoglu_simulation_2015}.

\subsubsection{Construction of a local interpolant}\label{sec: Interpolant_construction}
Inside each local support domain ${}_{l}\Omega$, an approximation function is defined as a sum of weighted shape functions
\begin{equation}
	{}_{l}y_{\xi}(\ve{p}) \approx \displaystyle\sum_{i=1}^{{}_{l}N} {}_{l}\alpha_{i, \xi} \: {}_{l}\Phi_{i} (\ve{p}),
\end{equation}
where index $i$ runs over ${}_{l}N$ nodes defined by ${}_{l}\Omega$ and $\xi$ runs over space dimensions of the approximated field $\xi = 1, ..., n_{d}$. Weight coefficients are represented by ${}_{l}\alpha_{i,\xi}$. For the shape functions, Radial Basis Functions (RBFs) are used. In general, RBFs are functions where the value of the function depends only on the distance from the RBF center ${}_{l}\ve{p}$. This can be written as ${}_{l}\Phi(\ve{p}) = \Phi(||\ve{p}-{}_{l}\ve{p}||)$. Depending on the class of RBFs the local interpolation problem may become ill-conditioned. This problem is overcome by adding monomials $p_i(\ve{p})$, $i=1, ..., M$ to the approximation \cite{fasshauer_meshfree_2007}. This combined approximation basis is also called augmented RBFs. 
\begin{equation}\label{eq: loc_app}
	{}_{l}y_{\xi}(\ve{p}) \approx \displaystyle\sum_{i=1}^{{}_{l}N} {}_{l}\alpha_{i,\xi} \: {}_{l}\Phi_{i}(\ve{p})  + \displaystyle\sum_{i=1}^{M} {}_{l}\alpha_{({}_{l}N+i), \xi} \:  p_{i}(\ve{p}) =
	\displaystyle\sum_{i=1}^{{}_{l}N+M} {}_{l}\alpha_{i, \xi} \: {}_{l}\Psi_{i}(\ve{p}).
\end{equation}
The above expression is written in a compact form where $\Psi_i(\ve{p})$, $i =1, ..., {}_{l}N +M$ represents the complete set of shape functions used, either RBFs ($i \le {}_{l}N$) or monomials ($i > {}_{l}N$). Local interpolation problem at each $l^\text{th}$ node can be written as a system of $n_{d}({}_{l}N + M)$ linear equations
\begin{equation}\label{eq: int_mat_problem}
	\displaystyle\sum_{\chi=1}^{n_{d}} \displaystyle\sum_{i=1}^{{}_{l}N + M} {}_{l}A_{ji,\xi \chi} \: {}_{l}\alpha_{i,\chi} = {}_{l}\gamma_{j,\xi},
\end{equation}
where the interpolation matrix ${}_{l}A_{ji,\xi \chi}$ is defined as 
\begin{equation}
	{}_{l}A_{ji, \xi \chi} = 
	\begin{cases*}
		\Psi_{i}({}_{l}\ve{p}_j) \delta_{\xi \chi} & if ${}_{l}\ve{p}_j \in \Omega \cup \Gamma$\\
		p_{j-{}_{l}N}({}_{l}\ve{p}_i) \delta_{\xi \chi} & if $j > {}_{l}N$ and $i \leq {}_{l}N$\\
		0         & otherwise
	\end{cases*},
\end{equation}
and the vector of known values as 
\begin{equation}
	{}_{l}\gamma_{j, \xi}=
	\begin{cases*}
		y_{\xi}({}_{l}\ve{p}_j) & if ${}_{l}\ve{p}_j \in \Omega \cup \Gamma$\\
		0         & otherwise.
	\end{cases*}
\end{equation}

\subsubsection{Approximation of differential operators}\label{sec: App_DO}
Application of linear differential operator $\mathcal{L}$ on local interpolation function, using the equation (\ref{eq: loc_app}), can be written as
\begin{equation}\label{eq: gen_op_act2}
	\mathcal{L} {}_{l}\bm{y}(\ve{p})_{\xi} = \displaystyle\sum_{\chi=1}^{n_d} \mathcal{L}_{\xi \chi} \: {}_{l}y_{\chi}(\ve{p}) \approx \displaystyle\sum_{\chi=1}^{n_d} \displaystyle\sum_{i=1}^{{}_{l}N+M} {}_{l}\alpha_{i, \chi} \: \mathcal{L}_{\xi \chi} \: {}_{l}\Psi_{i}(\ve{p}),
\end{equation}
where $\mathcal{L}$ preserves the order of the (generally tensorial) physical field, such as Laplacian operator $\nabla^2$  for example.
Since ${}_{l}\alpha_{i,\chi}$ are constants, the operator $\mathcal{L}$ acts only on shape functions. Constants from equation (\ref{eq: int_mat_problem}) can now be expressed as 
\begin{equation}
	{}_{l}\alpha_{i,\chi} = \displaystyle\sum_{\zeta=1}^{n_{d}} \displaystyle\sum_{j=1}^{{}_{l}N + M} {}_{l}A_{ij,\chi \zeta}^{-1} \: {}_{l}\gamma_{j,\zeta},
\end{equation}
and then used in equation (\ref{eq: gen_op_act2})
\begin{equation}
	\mathcal{L} {}_{l}\bm{y}(\ve{p})_{\xi} \approx \displaystyle\sum_{\zeta=1}^{n_{d}} \displaystyle\sum_{j=1}^{{}_{l}N + M} {}_{l}\gamma_{j,\zeta} \displaystyle\sum_{\chi=1}^{n_d} \displaystyle\sum_{i=1}^{{}_{l}N+M} {}_{l}A_{ij,\chi \zeta}^{-1} \: \mathcal{L}_{\xi \chi} \: {}_{l}\Psi_{i}(\ve{p}).
\end{equation}
A more compact form can be written as
\begin{equation} \label{eq: ADO}
	\mathcal{L} {}_{l}\bm{y}(\ve{p})_{\xi} \approx \displaystyle\sum_{\zeta=1}^{n_{d}} \displaystyle\sum_{j=1}^{{}_{l}N + M} {}_{l}\gamma_{j,\zeta} \: {}_{l}\mathscr{W}_{j,\xi \zeta}(\ve{p}),
\end{equation}
where ${}_{l}\mathscr{W}_{j,\xi \zeta}(\ve{p})$ represent operator coefficients, where position vector $\ve{p}$ is generally arbitrary.
Here, the similarity with the FD method can be seen, where the action of an operator on a function is evaluated as a weighted sum of known values and operator coefficients. These are determined in the preprocessing step (only once) for each ${}_{l}\Omega$ and differential operator.

\subsubsection{The choice of the radial basis function}\label{sec: PHS}
In this paper, a particular type of RBFs, called polyharmonic splines (PHSs) are used. In contrast with other commonly used RBFs, such as Gaussian or Multiquadric RBFs, with PHSs there is no need for complicated and time-consuming determination of an optimal shape parameter. It is well known \cite{fasshauer_meshfree_2007} that using only RBFs to construct an interpolant leads to stagnation error, meaning that the interpolation error does not decrease with node refinement. In \cite{flyer_role_2016}, it has been shown that the use of augmented PHSs decreases stagnation error and that the augmentation order controls the order of convergence. This makes PHSs very attractive to work with. The employed definition of PHS function, centered at the ${}_{l}\ve{p}_i$, where $i = 1, ... N$, that includes support size rescaling, is 
\begin{equation}
	{}_{l}\Phi_{i} (\ve{p}) = \left(\frac{||\ve{p}-{}_{l}\ve{p}_{i}||}{{}_{l}h}\right)^{m},
\end{equation}
where $m = 1, 3, 5, ...$ stands for the order of PHS, and ${}_{l}h$ represents an average distance from the central node, given as
\begin{equation}
	{}_{l}h = \sqrt{\displaystyle\sum_{i=2}^{{}_{l}N} \frac{||{}_{l}\ve{p}-{}_{l}\ve{p}_{i}||^2}{{}_{l}N - 1 }}.
\end{equation}
Some important observations from \cite{bayona_role_2017} should be emphasized here. In order to obtain the $p^\text{th}$ order of convergence with respect to $h$, where $p$ represents the maximum order of monomials used in augmentation and when it holds that $m < p$, the condition ${}_{l}N \gtrsim 2 M$, where $M = \binom{p+n_d}{p}$ should be satisfied. This means that the number of support nodes should be at least two times larger than needed for the approximation with $p^\text{th}$-order monomials.

\subsection{Incremental solution of boundary value problem}
Stress state is, in elasto-plasticity, a function of strain state and a history of strains to which the body has been subjected. Due to that, the elasto-plastic BVP is path-dependent and should be solved incrementally, where the external load is applied incrementally. Suppose that in the increment $n$, the internal stresses, which are appropriately determined through constitutive relations at each material point in the domain, are in equilibrium with the external load. Next, a new load increment $n+1$ is applied, and the equilibrium condition is no longer satisfied. The task of analysis is now to find the correct internal force, that is in equilibrium with current external load. This can be for some material point written as
\begin{equation}\label{eq: ravnovesje}
	\ve{r} (\ve{u}_{n+1}) = \ve{f}^{int} |_{n+1} - \ve{f}^{ext} |_{n+1} = 0,
\end{equation}
where $\ve{f}^{ext}$ represents an external load and $\ve{r}$ is a residual that should be equal to zero at the end of the load increment. Internal force $\ve{f}^{int}$ can be, using the equation (\ref{eq: balance DE}), expressed as a function of stress tensor as
\begin{equation}\label{eq: Internal_force}
	\ve{f}^{int} = \nabla \cdot \bm\sigma(\nabla^s \ve{u}).
\end{equation}
Due to the non-linear relationship between stresses and strains, the residual is linearized by Taylor expansion. The linearized form can be written as
\begin{equation} \label{eq: N-R1}
	\frac{\partial }{\partial \ve{u}} \big[\nabla \cdot(\bm\sigma(\nabla^s \ve{u})) \big] \bigr\rvert_{n+1}^{(i-1)} \delta \ve{u}^{(i)} =  -\ve{r}|_{n+1}^{(i-1)} , 	
\end{equation}
using the first-order expression for the internal force. It is solved for $\delta \ve{u}^{(i)}$. Index $n$ is referred to the loading step and $i$ to the iteration index of the iterative solving of a non-linear equation. After each iteration, the increment of the displacement vector is updated as $\Delta \ve{u}^{(i)}  = \Delta \ve{u}^{(i-1)} + \delta \ve{u}^{(i)}$, where the initial guess is $\Delta \ve{u}^{(0)} = 0$. This algorithm, also called NRIA, keeps iterating until the residual is sufficiently small. Then the displacement vector is updated as $\ve{u}_{n+1} = \ve{u}_{n} + \Delta \ve{u}$. To solve equation (\ref{eq: N-R1}), the Jacobian on the left-hand side must be known. Taking the derivative over displacement, it holds that
\begin{equation}\label{eq: TTM0}
	\frac{\partial }{\partial \ve{u}} \big[\nabla \cdot(\bm\sigma(\nabla^s \ve{u})) \big] \bigr\rvert_{n+1}^{(i-1)} = \nabla \cdot \left(\frac{\partial \bm\sigma}{\partial \bm\varepsilon} \nabla^s \right) \biggr\rvert_{n+1}^{(i-1)} = \nabla \cdot \left(\textbf{\textsf{D}} \nabla^s \right) \bigr\rvert_{n+1}^{(i-1)},
\end{equation}
where $\textbf{\textsf{D}}=\partial \bm\sigma / \partial \bm\varepsilon$ is the tangent operator, which must be consistent with the integration scheme used for solving constitutive equations (\ref{eq: Hooks Law}-\ref{eq: K-T}) to ensure the convergence of the NRIA \cite{simo_computational_2000}. It is called a consistent tangent operator and is for implicit integration of constitutive equations, used in this work elaborated in Section \ref{sec: RMA}. The final equation of interest is obtained as
\begin{equation}\label{eq: MAIN_NR}
	\nabla \cdot \left(\textbf{\textsf{D}} \nabla^s \right) \bigr\rvert_{n+1}^{(i-1)} \delta \ve{u}^{(i)} = -\ve{r} |_{n+1}^{(i-1)}. 
\end{equation}
by combining equations \ref{eq: N-R1} and \ref{eq: TTM0}. Discretized form of differential operators is defined on the reference configuration and is not changing during increments, as usual in small strain approximation. With discretized equation \ref{eq: MAIN_NR} at hand, a global system of equations can be composed.
By obtaining the solution of $\Delta \ve{u}^{(i)}$ at each node of the observed domain, the strain tensor increment $\Delta \bm\varepsilon|_{n+1}^{(i)}$ can be determined by equation (\ref{eq: TSDeformacij}). This serves as an input for the integration model (presented in Section \ref{sec: RMA}), where the constitutive equations are solved locally at each node. From there, a new stress state $\bm\sigma |_{n+1}^{(i)}$ and other state variables are calculated. Finally, with the stress state calculated, a new residual is determined as $\ve{r} |_{n+1}^{(i)} = \nabla \cdot \bm\sigma |_{n+1}^{(i)} - \ve{f}^{ext} |_{n+1}$ and checked as $||\ve{r} |_{n+1}^{(i)}||/||\ve{f}^{ext} |_{n+1}|| \leq e_{tol}^{NR}$, where $e_{tol}^{NR}$ represents equilibrium convergence tolerance. If the condition is not satisfied NRIA index is updated $i=i+1$, and a new cycle is performed.

\subsubsection{Return Mapping Algorithm - RMA}\label{sec: RMA}
\begin{figure}[H]
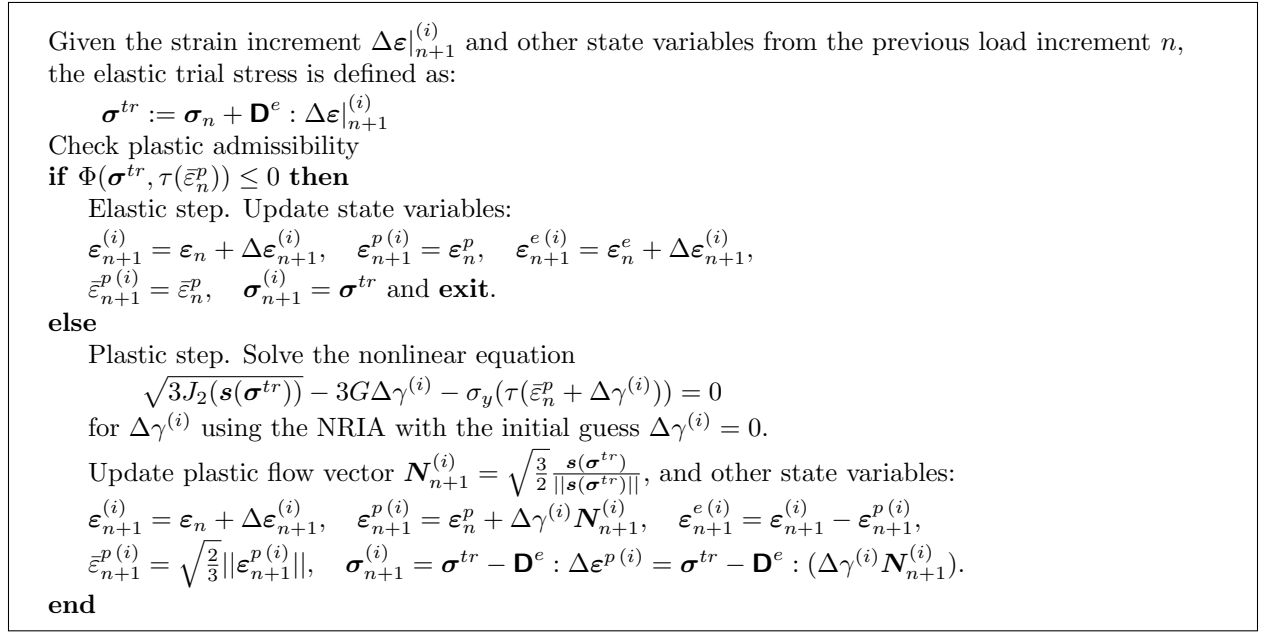

	\begin{algorithm}[H]
		\DontPrintSemicolon
		\BlankLine
		Given the strain increment $\Delta \bm\varepsilon|_{n+1}^{(i)}$ and other state variables from the previous load increment $n$, \; the elastic trial stress is defined as: \;
		$\qquad \bm\sigma^{tr} := \bm\sigma_{n} + \textbf{\textsf{D}}^e : \Delta \bm\varepsilon|_{n+1}^{(i)}$\;
		Check plastic admissibility\;
		\eIf{$\Phi(\bm\sigma^{tr}, \tau(\bar{\varepsilon}^p_{n})) \leq 0$}
		{
			Elastic step. Update state variables:\;
			$\bm\varepsilon_{n+1}^{(i)} = \bm\varepsilon_n + \Delta \bm\varepsilon_{n+1}^{(i)}, \quad \bm\varepsilon_{n+1}^{p\,(i)} = \bm\varepsilon_n^p , \quad \bm\varepsilon_{n+1}^{e \, (i)} = \bm\varepsilon_{n}^e + \Delta \bm\varepsilon_{n+1}^{(i)},$\;
			$\bar{\varepsilon}^{p \, (i)}_{n+1} = \bar{\varepsilon}^p_{n}, \quad \bm\sigma_{n+1}^{(i)} = \bm\sigma^{tr}$ and \textbf{exit}.
		}{
			Plastic step. Solve the nonlinear equation\; 
			$\qquad \sqrt{3 J_2 (\bm{s}(\bm\sigma^{tr}))} - 3 G \Delta \gamma^{(i)} - \sigma_y(\tau(\bar{\varepsilon}^p_{n} + \Delta \gamma^{(i)})) = 0$\;
			for $\Delta \gamma^{(i)}$ using the NRIA with the initial guess $\Delta \gamma^{(i)} = 0$. \;
			Update plastic flow vector $\bm{N}_{n+1}^{(i)} = \sqrt{\frac{3}{2}} \frac{\bm{s}(\bm\sigma^{tr})}{||\bm{s}(\bm\sigma^{tr})||}$, and other state variables:\;
			$\bm\varepsilon_{n+1}^{(i)} = \bm\varepsilon_n + \Delta \bm\varepsilon^{(i)}_{n+1}, \quad \bm\varepsilon_{n+1}^{p \, (i)} = \bm\varepsilon_n^p + \Delta \gamma^{(i)} \bm{N}_{n+1}^{(i)}, \quad \bm\varepsilon_{n+1}^{e \, (i)} = \bm\varepsilon_{n+1}^{(i)} - \bm\varepsilon_{n+1}^{p \, (i)},$\;
			$\bar{\varepsilon}^{p \, (i)}_{n+1} = \sqrt{\frac{2}{3}}||\bm\varepsilon_{n+1}^{p \, (i)}||, \quad \bm\sigma_{n+1}^{(i)} = \bm\sigma^{tr} - \textbf{\textsf{D}}^e : \Delta\bm\varepsilon^{p \, (i)} = \bm\sigma^{tr} - \textbf{\textsf{D}}^e : (\Delta \gamma^{(i)} \bm{N}_{n+1}^{(i)})$.
		}
		
	\end{algorithm}
	\caption{RMA for the von Mises associate model with non-linear isotropic hardening.}
	\label{alg1}
\end{figure}

In this section a brief revision of the algorithm for solving a set of constitutive equations (\ref{eq: Hooks Law} - \ref{eq: K-T}) is given. An appropriate integration algorithm must be employed to solve them. In this work the two step elastic predictor/plastic corrector algorithm is used. The steps are as follows: in the first step, given a strain increment, the stress state, which is a function of accumulated plastic strain, is predicted as purely elastic. Then, if conditions in equation (\ref{eq: K-T}) are not satisfied, a set of equations (\ref{eq: additive_split} -\ref{eq: K-T}) is iteratively solved with RMA to correct state variables properly. RMA is implemented with an implicit Euler scheme that is unconditionally stable.
The procedure is summarized in Figure \ref{alg1}. For von Mises associated isotropic hardening, a non-linear scalar equation must be solved for $\Delta \gamma$, where all non-linearity comes from the hardening curve. If the hardening curve is linear or constant, state variables in the plastic step can be determined exactly.

For this case, the consistent tangent operator is derived as  
\begin{equation}\label{eq: CTO}
	\textbf{\textsf{D}} \bigr\rvert_{n+1}^{(i)}=
	\begin{cases*}
		\textbf{\textsf{D}}^e         & if $\Delta \gamma^{(i)} = 0$\\
		\textbf{\textsf{D}}^e - \frac{6 G^2 \Delta \gamma^{(i)}}{\sqrt{3 J_2 (\bm{s}(\bm\sigma^{tr}))}} \textbf{\textsf{I}}_d + 6 G^2 \left(\frac{\Delta \gamma^{(i)}}{\sqrt{3 J_2 (\bm{s}(\bm\sigma^{tr}))}} - \frac{1}{3 G + H} \right) \bar{\bm{N}}_{n+1}^{(i)} \otimes \bar{\bm{N}}_{n+1}^{(i)} & if $\Delta \gamma^{(i)} > 0$
	\end{cases*} ,
\end{equation}
where $H = \partial \sigma_y(\tau(\bar{\varepsilon}_n^p + \Delta \gamma))/\partial \bar{\varepsilon}^p$ represents the slope of the hardening curve, also called the hardening modulus. $\bar{\bm{N}} = \bm{N}/||\bm{N}||$ stands for unit plastic flow direction and $\textbf{\textsf{I}}_d$ for deviatoric projection tensor \cite{de_souza_neto_computational_2011}.

\subsection{Solution approaches}
The function is expected to be continuous when applying the differential operator to it with RBF-FD, as shown in equation \ref{eq: ADO}. Observing equation \ref{eq: MAIN_NR}, the divergence operator acts on the stress field, where $\textbf{\textsf{D}}$ (see equation \ref{eq: CTO}) is generally discontinuously changing in space because of the discontinuous transition from elastic to plastic region, and vice versa. This discontinuity represents an essential challenge when solving the elasto-plastic BVP with RBF-FD.

This work presents three different approaches of discretizing equation \ref{eq: MAIN_NR} with the RBF-FD method. Discretization of BCs is, in all cases, the same. All approaches are appropriate and successful in resolving a purely elastic response. However, it will be shown in the following sections that the \textit{direct} and \textit{composed} approaches cannot be used to overcome stability and accuracy problems when solving elasto-plasticity.

\subsubsection{Direct approach}
In this approach, the divergence operator in equation \ref{eq: MAIN_NR} is evaluated analytically. This results in the following expression
\begin{equation}\label{eq: Direct Approach}
	\left[ \left(\nabla \cdot \textbf{\textsf{D}}\right) : \nabla^s + \textbf{\textsf{D}} : \nabla \otimes \nabla^s \right]  \bigr\rvert_{n+1}^{(i-1)}   \delta \ve{u} = -\ve{r}|_{n+1}^{(i-1)}.
\end{equation}
Here local discretization is performed by calculating operator coefficients of $\nabla \cdot$, $\nabla^s$ and $\nabla \otimes \nabla^s$ directly via RBF-FD. According to equation \ref{eq: Direct Approach}, the coefficients are computed for each node and then expanded into the global stiffness matrix. The internal force is computed as defined by equation \ref{eq: Internal_force}, for which divergence operator coefficients must be provided. In this approach, $2^{\text{nd}}$-order derivatives must be discretized. 


\subsubsection{Composed approach}
In the \textit{composed} approach, the divergence operator in equation \ref{eq: MAIN_NR}  is evaluated in the discretized form. Here the discrete stress field is determined as
\begin{equation}\label{eq: add_s_field}
	\delta \bm\sigma |_{n+1}^{(i)} = (\textbf{\textsf{D}} \nabla^s) |_{n+1}^{(i-1)} \delta \ve{u}^{(i)},
\end{equation}
where the discretized form of the term $(\textbf{\textsf{D}} \nabla^s) |_{n+1}^{(i-1)}$ is afterwards globally arranged in a rectangular sparse matrix $\bm{K}_{\sigma}$. For $\delta \bm\sigma$, the equilibrium equation must be satisfied
\begin{equation}
	\nabla \cdot \delta \bm\sigma |_{n+1}^{(i)} = -\ve{r} |_{n+1}^{(i-1)},
\end{equation}
where the discretized form of the divergence operator is composed into another rectangular sparse matrix $\bm{K}_{div}$. The final global tangent stiffness matrix is then composed as a product $\bm{K}_T = \bm{K}_{div} \, \bm{K}_{\sigma}$. The internal force vector is calculated in the same way as in the \textit{direct} approach. Unlike in the \textit{direct} approach, only first-order derivatives must be discretized here.


\subsubsection{Hybrid approach}\label{sec: HybridApp}
\begin{figure}
	\begin{centering}
		\includegraphics[scale=0.9]{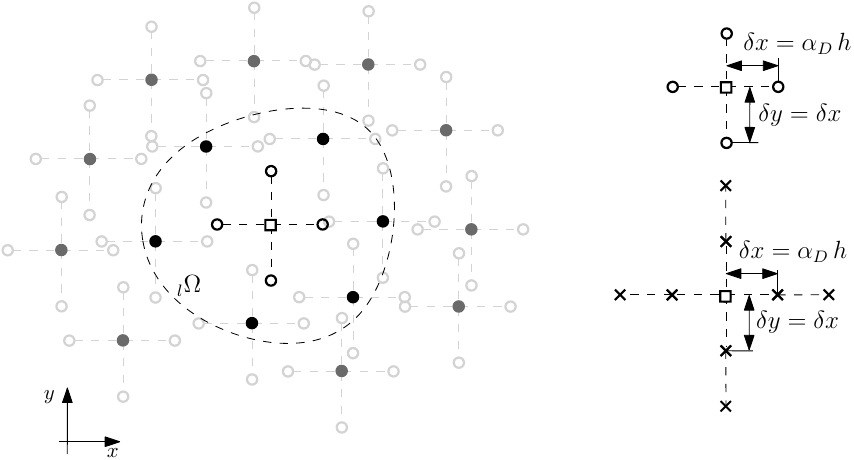}
		\caption{Left: Arrangement of collocation nodes with $2^{\text{nd}}$-order FD stencils assigned. The collocation node, denoted with a square, represents the central node of the local support domain ${}_{l}\Omega$. Right: $2^{\text{nd}}$-order (top) and $4^{\text{th}}$-order (bottom) FD stencils with prescribed dimensions.}
		\label{fig: new_staff_1}
	\end{centering}
\end{figure}

The goal of the \textit{hybrid} approach is to evaluate the divergence operator with the FD method. A virtual FD stencil is prescribed to each collocation node, as shown in Figure \ref{fig: new_staff_1} (left). At the FD stencil nodes, called here secondary nodes (SNs), all secondary variables (stresses and strains) are evaluated via RBF-FD. As in the \textit{composed} approach, coefficients for the term $(\textbf{\textsf{D}} \nabla^s) |_{n+1}^{(i-1)}$ are determined first. For the case of $p_{FD}=2$, where $p_{FD}$ represents the order of the FD method, four different sets of coefficients are determined within each ${}_{l}\Omega$, each set belonging to a different secondary node. Then, the divergence in the collocation node is expressed by manipulating of the determined coefficients in a FD manner. With coefficients computed, the global stiffness matrix is composed as in the $\textit{direct}$ approach. Since the stress values are computed on secondary nodes, the internal force vector is also computed via the FD method. It should be emphasized that when determining the coefficients of $(\textbf{\textsf{D}} \nabla^s) |_{n+1}^{(i)}$, the local support domain for an individual secondary node does not change but remains the same as the one belonging to the central collocation node.

With the introduction of the FD stencil, the question of its size arises. As shown in Figure \ref{fig: new_staff_1} (right), the size of the FD stencil is defined as a fraction of the distance between the two closest collocation nodes. For $p_{FD} = 2$ it is defined as $2\delta x= 2\delta y = 2 \alpha_D h$ and for $p_{FD} = 4$ as $4\delta x= 4\delta y = 4 \alpha_D h$. The value of the newly introduced parameter $\alpha_D$ is studied in the following sections. Since the FD method is used, the coefficients for determining the divergence operator are already known, which reduces computational complexity, as there is no need to determine the coefficients with the RBF-FD method. What increases the computational complexity compared to the \textit{direct} and composite methods is solving the material iteration in all SNs. A detailed description of the discretization with a \textit{hybrid} approach is given in \ref{app: Operator_der}. 

\subsection{Stabilization of boundary conditions}\label{sec: Stabliization BC}
It is well known that one-sided stencils at (and near) the boundaries lead to the loss in accuracy and stability induced by spurious oscillations \cite{fasshauer_meshfree_2007}. This problem arises especially when derivatives on the boundary, as in the case of the Neumann BC (at $\Gamma_F$ and $\Gamma_T$), are specified. Many different approaches \cite{bayona_role_2017, gerace_model-integrated_2014, shu_local_2003, zheng_meshfree_2015, zheng_local_2018, slak_generation_2019, khosrowpour_strong-form_2019} were proposed to ensure a more accurate and stable evaluation of normal derivatives ($\partial/\partial \ve{n}$), but it still remains a challenge to this day. In one of the first approaches of computing the derivative on the boundary, a one-sided second-order FD scheme was used \cite{shu_local_2003, gerace_model-integrated_2014}. The first two layers of nodes inside the domain were positioned on the opposite side of the outward-facing normal vector. A similar approach was proposed with RBF-FD approximation along a single line, where more than two layers were added, which is difficult to achieve when dealing with complex geometries. However, when possible, it provides a higher order of accuracy \cite{zheng_meshfree_2015}. In \cite{zheng_meshfree_2015}, two additional approaches were proposed. In the first, Neumann BCs are evaluated as a combination of two RBF-FD approximations along two different lines, which also suffers from inflexibility. The second approach uses a set of fictitious nodes that are placed along the normal vector inside the domain. In this case, the values are being interpolated from CNs and used in the same single-line RBF-FD manner. The last approach is the most flexible but also most unstable because of the interpolation on fictitious nodes \cite{zheng_meshfree_2015}. A similar approach was introduced where only one fictitious node (one layer) is generated. Then the normal derivative is calculated with a one-sided first-order FD stencil where both fictitious and boundary nodes are used \cite{zheng_local_2018}. Another popular and well-established way arising from the FD method is the use of ghost nodes, where additional layer(s) of nodes are added outside the domain in the direction of outward-facing normal vectors. This allows higher-order central FD schemes on the boundary, which increases stability but comes with the cost of introducing additional nodes. When using RBF-FD, ghost nodes allow more accurate approximation and higher stability since each boundary node lies closer to the support domain center \cite{bayona_role_2017, slak_generation_2019, khosrowpour_strong-form_2019}.

In \cite{bayona_role_2017}, augmented PHSs for solving elliptic PDE were used. It was found that the accuracy and stability problems can be overcome with a sufficient increase in the stencil size. Nevertheless, when solving non-linear elliptic PDE, the ghost nodes were used to ensure the solvability of the investigated problem.  

\begin{figure}
	\begin{centering}
		\includegraphics[scale=1.1]{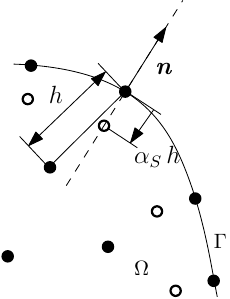}
		\caption{Shifted boundary nodes (hollow circles) orthogonally to the boundary to evaluate discretized differential operators.}
		\label{fig: new_staff_2}
	\end{centering}
\end{figure}

In this work, augmented PHSs are used. In contrast to \cite{bayona_role_2017}, the authors of this work found that stabilization of BCs with ghost nodes has no significant effect in this particular case. Instead, a similar approach to the one presented in \cite{zheng_local_2018} is employed here. As seen in Figure \ref{fig: new_staff_2}, BCs are not evaluated on the boundary nodes but on the virtually shifted boundary nodes presented with hollow circles. These are obtained by virtually shifting boundary nodes in the opposite direction of the outward-facing normal vectors for the distance of $\alpha_S \, h$, where $\alpha_S$ is a non-negative scalar value. With this approach, differential operators are evaluated on of fictitious nodes that do not lay exactly on $\Gamma$. This approach is then applied to each boundary section where any kind of derivative is needed to evaluate the BCs. The disadvantage of this approach is that the BC is not satisfied precisely on the boundary, but rather on the fictitious node.

As proposed in \cite{bayona_role_2017}, an increase in stencil size only on $\Gamma$ should also work. This seems like a good option, but only under the assumption that the solution field is smooth. As previously mentioned, this is not the case here since the consistent tangent operator is discontinuous. 

\subsection{Numerical implementation}
The numerical implementation of the presented approaches is written in the Fortran 2018 programming language and compiled with Intel Fortran Compilers Classic 2021.1.1. Computation was performed with an Intel(R) Core(TM) i7-8750H CPU containing six cores with a base clock speed of 2.20 GHz. 

\section{Numerical performance in the elastic range} \label{sec: NUM_EX_EL}
The performance of the proposed approaches is first investigated within a linear elastic range. Three benchmarks with different levels of complexity are studied in parallel. First, a short comment is given on a well-established (within linear elastic range) \textit{direct} approach. Then, the performance of the \textit{composed} approach is studied in terms of convergence. In the following, the main focus is given to the \textit{hybrid} approach. The effect of the FD stencil size $\alpha_D$ is investigated regarding convergence, stability and accuracy. At the end of the section, the effect of the BC stabilization parameter $\alpha_S$ is studied in terms of stability, accuracy and convergence. 



\subsection{Elastic benchmarks}

\subsubsection{Timoshenko beam}
The case often used to verify the implementation of methods in solid mechanics is the Timoshenko beam \cite{atluri_meshless_2004, liu_meshfree_2009, simonenko_optimal_2014, mavric_meshless_2017}. The solution of the displacement field is a polynomial of a $3^{\text{rd}}$ order with mixed terms. For a state of plane stress, it is defined in Cartesian coordinates ($x_1, x_2$) as 
\begin{equation}
	\begin{split}
		&\hat{u}_1(x_1, x_2) = \frac{F_0}{6EI}\left(x_2 - \frac{D}{2}\right) \left((6L - 3 x_1)x_1 + (2 + \nu) \left(x_2^2 - D x_2\right) \right),\\
		&\hat{u}_2(x_1, x_2) = -\frac{F_0}{6EI}\left(3 \nu \left( x_2 - \frac{D}{2} \right)^2 \left( L-x_1 \right)  + \left( 4+5\nu \right) \frac{D^2 x_1}{4} + (3L -x_1)x_1^2 \right), 
	\end{split} \label{eq: Timo}
\end{equation}
where $I = D^3/12$ represents the $2^{\text{nd}}$ moment of area, $D$ is the height, and $L$ is the length of a beam \cite{s_timoshenko_theory_1951}. Values of parameters used and geometry with specified BCs are shown for all elastic benchmarks in Figure \ref{fig: EL_ALL}.

\subsubsection{Stretching of a plate with a circular hole}\label{sec: PWH}
The next case, also commonly used to verify the numerical method \cite{atluri_meshless_2004, liu_meshfree_2009, mavric_meshless_2017}, is a section of the infinite plate with a circular hole subjected to a constant load $\sigma_{\infty}$ as $x_1 \to \infty$.
The exact solution in polar coordinates $(r, \theta)$ is expressed as
\begin{equation} 
	\begin{split}
		&\hat{u}_1(r, \theta) = \frac{\sigma_{\infty} \, R_i}{8G} \left( \left( \beta+1 \right) \frac{r}{R_i}\text{cos}(\theta) + \frac{2 R_i}{r} \left(  \left(\beta+1\right) \text{cos}(\theta) + \text{cos}(3\theta) \right) -2\left(\frac{R_i}{r}\right)^3 \text{cos}(3\theta) \right), \\
		&\hat{u}_2(r, \theta) = \frac{\sigma_{\infty} \, R_i}{8G} \left( \left( \beta-3 \right) \frac{r}{R_i}\text{sin}(\theta) + \frac{2 R_i}{r} \left(  \left(1-\beta\right) \text{sin}(\theta) + \text{sin}(3\theta) \right) -2\left(\frac{R_i}{r}\right)^3 \text{cos}(3\theta) \right),
	\end{split} \label{eq: PLW_sol}
\end{equation}
where $\beta = (3-\nu)/(1+\nu)$ for a plane stress case and $R_i$ denotes the radius of a hole. Since the solution is symmetric over the $x_1$ axis and anti-symmetric over the $x_2$ axis, only one-quarter of the plate is considered. Traction load is obtained from analytical stress field solution, derived from equation \ref{eq: PLW_sol}.


\subsubsection{Internally pressurized annulus} \label{sec: IPA}
The last elastic case considered is an internally pressurized annulus \cite{s_timoshenko_theory_1951}. The exact solution of radial displacement for a plane strain state is defined as 
\begin{equation}
	\hat{u}_r(r) = -\frac{R_i^2 p_0}{ 2(G + \lambda)\left( R_i^2 - R_o^2 \right)  } r  -\frac{R_i^2 R_o^2 p_0}{ 2 G \left( R_i^2 - R_o^2 \right) } \frac{1}{r},
\end{equation}
where $R_i$,  $R_o$ and $p_0$ present the inner radius, outer radius, and inner pressure, respectively. Because the analytical solution has cylindrical symmetry, we solved the problem on a one-quarter of the annulus.

\begin{figure}
	\begin{centering}
		\includegraphics[width=1\textwidth]{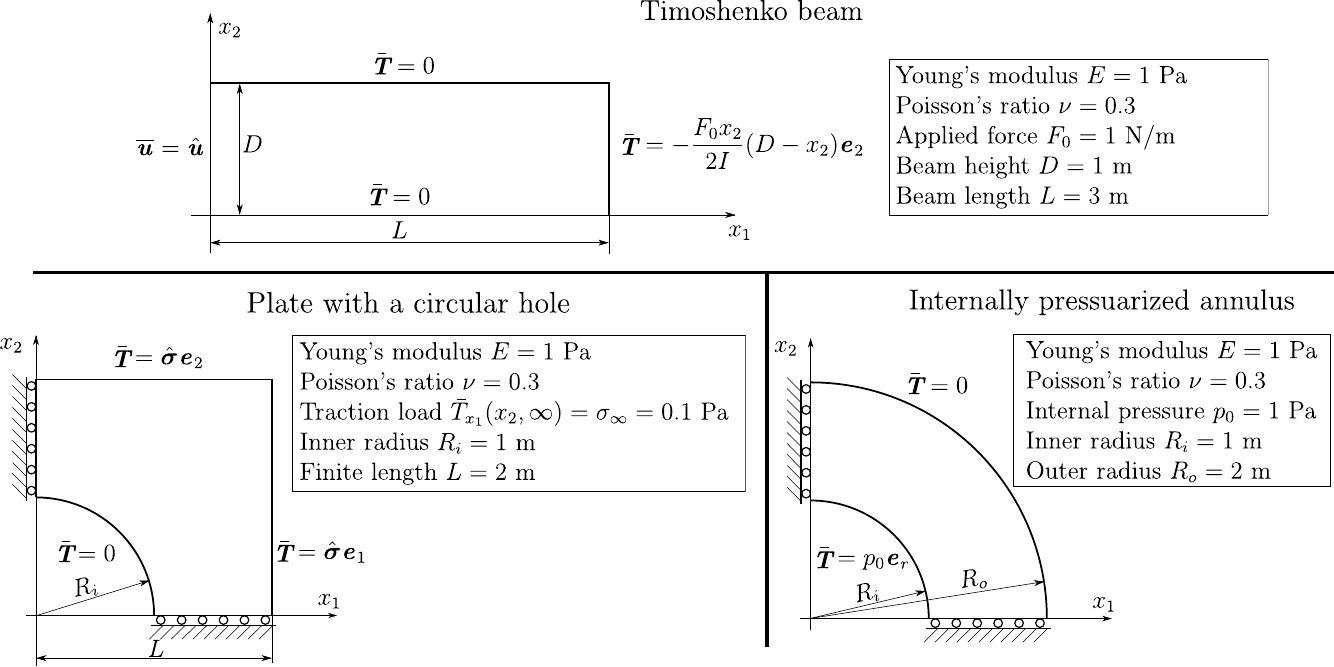} 
		\caption{Schemes of elastic benchmarks with geometry, boundary conditions, and material parameters used in the calculation. Top: Timoshenko beam, bottom-left: plate with a circular hole, bottom-right: internally pressurized annulus.}
		\label{fig: EL_ALL}
	\end{centering}
\end{figure}

\subsection{Direct approach performance} \label{sec: Direct performance}
The \textit{direct} approach has already been successfully verified using PHSs in works such as \cite{slak_adaptive_2020} and \cite{slak_generation_2019}. They have shown, as initially proposed in \cite{flyer_role_2016, bayona_role_2017}, that the order of monomial augmentation governs the convergence order.

\subsection{Composed approach performance} \label{sec: Composed performance}
Similarly, as in \cite{slak_adaptive_2020}, the convergence of the \textit{composed} approach is studied here. In this study, only SNAs are used since this is one of the main advantages of the proposed method. Examples of SNAs of the benchmarks are shown in \ref{app: Meshes}. The present study is performed for three different orders of augmentation $p = 2,3,4$, where the number of the support nodes is determined by ${}_{l}N = 2 M + 1$, since the proposed relation in \cite{flyer_role_2016} is ${}_{l}N \gtrsim 2 M$. The order of PHSs used is $m=3$. Under refinement, RBFs play an insignificant role in the accuracy of the approximation, so the order of PHSs matters little \cite{flyer_role_2016}. As mentioned in Section \ref{sec: PHS}, it should hold that $m<p$. When choosing $m>3$, larger stencils should be composed, which is not desirable from the computational complexity. For example, if $m=5$, this leads to $p=6 \rightarrow M=29 \rightarrow {}_{l}N \gtrsim 58$.  Here with $m=3$ and $p=2$ or $p=3$, the condition $m<p$ is not fulfilled. Nevertheless, the authors still tend to observe that regime of parameter use since low orders of augmentation induce smaller local interpolation problems significantly impact computational complexity at large (number of CNs) problems. No stabilization techniques are used to stabilize BCs, and no conditions on a maximum number of boundary nodes inside support domains are imposed. The employed error estimate is the relative $L_2$ norm $e_2$ defined by the following equation 
\begin{equation}\label{eq: e2}  
	e_{2} = \sqrt{\frac{\sum_{l=1}^{N} \left\| {}_{l}\ve{u} - {}_{l}\hat{\ve{u}} \right\|^2 }{\sum_{l=1}^{N} \left\| {}_{l}\hat{\ve{u}} \right\|^2}},
\end{equation}
where $\hat{\ve{u}}$ represents the exact and $\ve{u}$ the numerically obtained displacement value. For comparison with a well-established FEM, a plate with a circular hole case is studied with the Abaqus program package \cite{smith_abaqusstandard_2009}.  6-node quadratic finite elements (FEs) are employed with the average length of an FE's side adjusted to achieve a comparable number of unknowns as in RBF-FD cases with a consistent average node spacing.

\begin{figure}
	\begin{centering}
		\includegraphics[width=1\textwidth]{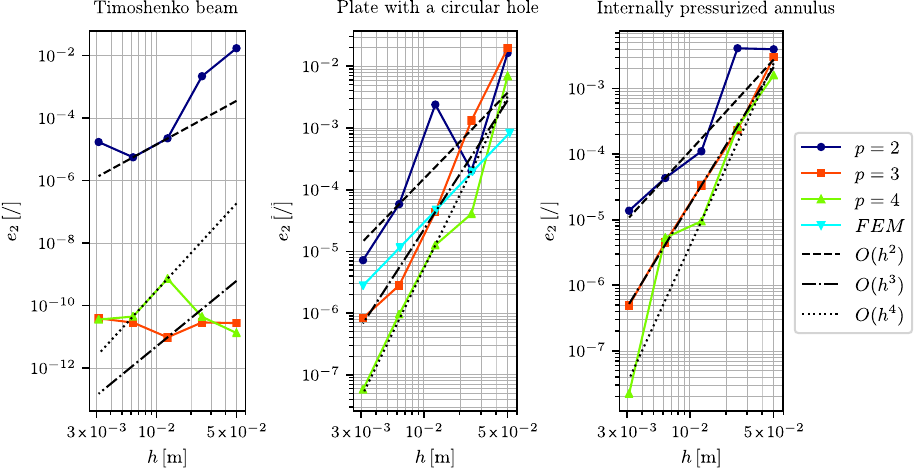} 
		\caption{Convergence of \textit{composed} approach for Timoshenko beam (left),  plate with a circular hole (center) and internally pressurized annulus (right).}
		\label{fig: Convergence_composed}
	\end{centering}
\end{figure}

Figure \ref{fig: Convergence_composed} presents the results of the convergence study on proposed elastic benchmarks using the \textit{composed} discretization approach. Here relative error is plotted as a function of the average node spacing. Straight lines present values of $O(h^p)$. As in the \textit{direct} approach, one can see that the augmentation order governs the convergence order. An exception can be seen in the Timoshenko beam when using the $3^{\text{rd}}$ and $4^{\text{th}}$ order of augmentation. No convergence in error is observed since the exact solution (polynomial of the $3^{\text{rd}}$ order) can be accurately determined using $3^{\text{rd}}$-order augmented basis functions. The presence of jumps that appear when increasing (or decreasing) discretization density is a well known phenomenon in RBF-FD field and the consequence of different SNAs. Many calculations should be performed with different SNAs to determine the range of possible solutions. An example of a procedure can be found in \cite{slak_adaptive_2020}. Here this is not done since a single calculation can demonstrate convergent behavior, which is the main point of interest here. 

Results show that FEM solution is converging exactly with the $2^{\text{nd}}$ order and is, compared to composed results with $p=2$, by nearly a factor of $\sim 10$ more accurate. Compared to solutions with $p=3, 4$, FEM solution is more accurate at small $h$, but composed solutions eventually get more accurate with refinement since they possess higher order of convergence.

\subsection{Hybrid approach performance}
In this section, the performance of the \textit{hybrid} approach is investigated. With the newly introduced parameter $\alpha_D$, specifying the size of a FD stencil, the convergence study for different values of $\alpha_D$ is presented first. The study is performed for different orders of RBF augmentation $p$ and different FD orders of divergence operator discretization $p_{FD}$. Then, the effect of $\alpha_D$ is studied in terms of stability and accuracy at different orders of FD method $p_{FD} = \left\{2,4\right\}$.


\subsubsection{Convergence of hybrid approach}\label{sec: Con_HY}
\begin{figure}[H]
	\begin{centering}
		\includegraphics[width=1\textwidth]{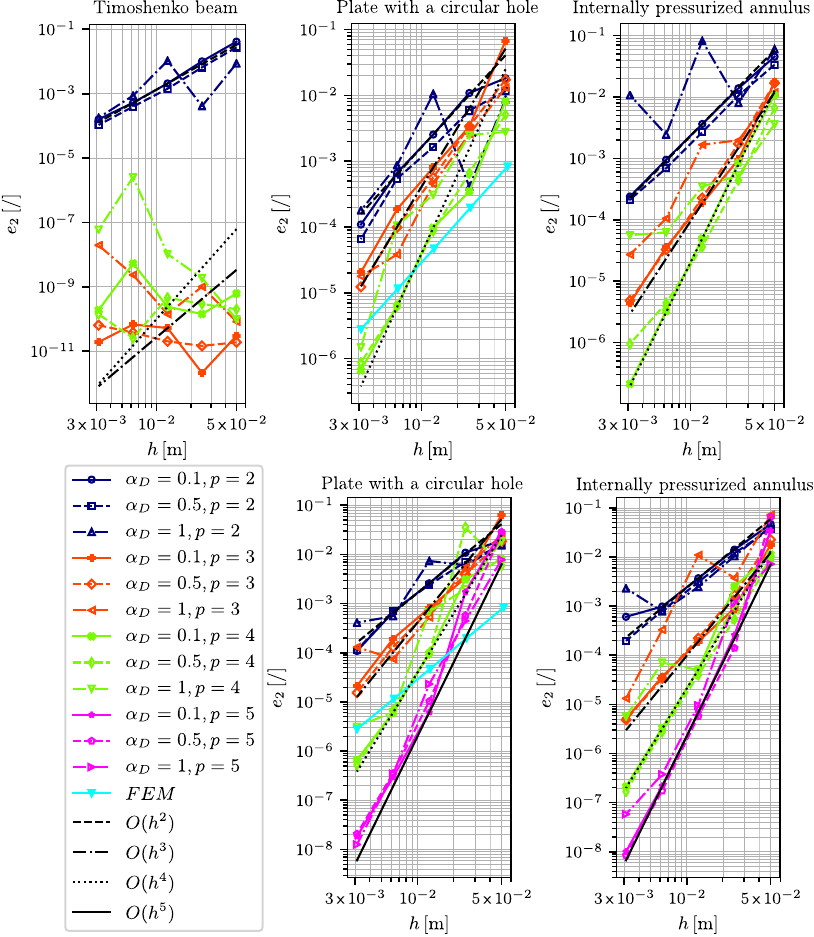} 
		\caption{Convergence of \textit{hybrid} approach using $2^{\text{nd}}$ order (top row) and $4^{\text{th}}$ order (bottom row) of an FD stencil. Timoshenko beam (left), plate with a circular hole (center) and internally pressurized annulus (right).}
		\label{fig: hybrid_ALL}
	\end{centering}
\end{figure}

The convergence study is performed for three different values of $\alpha_D = \left\{0.1, 0.5, 1\right\}$. 
In Figure \ref{fig: hybrid_ALL} (top row), the results are presented for the case of $p_{FD} = 2$. One can see that for $\alpha_D=0.1$ and $\alpha_D=0.5$, the convergence order is governed by the RBF augmentation. Interestingly, the jumps in the \textit{direct} and \textit{composed} approach mostly disappear with $\alpha_D=0.1, 0.5$. In that sense, the method becomes more stable in terms of convergence. On the other hand, effectively comparing error values with the \textit{composed} approach, presented in Figure \ref{fig: Convergence_composed}, the error in the \textit{hybrid} approach is approximately ten times larger. Compared with FEM results the error is almost 100 times larger. With $\alpha_D = 1$, the convergence order is spoiled, more jumps are present, and in the Timoshenko beam ($p=3, 4$), the error even starts increasing, where it is expected to experience no noticeable change. In Figure \ref{fig: hybrid_ALL} (bottom row), results for $p_{FD}=4$ are presented where a case with $p=5$ is added. Due to its relatively simple solution, the Timoshenko beam is not considered here. With small values of $\alpha_D$ again, the RBF augmentation order governs the order of convergence.

To sum up, to observe the $p^{\text{th}}$ -order of convergence, $\alpha_D < 1$ should be chosen.

%

%

\subsubsection{Size and discretization order of the FD stencil} \label{sec: A_D}
To investigate the impact of the FD stencil size and order on the proposed approach, the condition number $\kappa$ of the tangent stiffness matrix $\bm{K}_T$ and relative error $e_2$ are calculated for different values of $\alpha_D$. Condition number represents sensitivity to numerical noise and is defined as $\kappa = \sigma_{max}(\bm{K}_T) / \sigma_{min}(\bm{K}_T)$, where $\sigma_{max}$ and $\sigma_{min}$ are maximum and minimum singular values \cite{trefethen_numerical_1997}. The relative error is defined by equation \ref{eq: e2}.

\begin{figure}[H]
	\begin{centering}
		\includegraphics[width=1\textwidth]{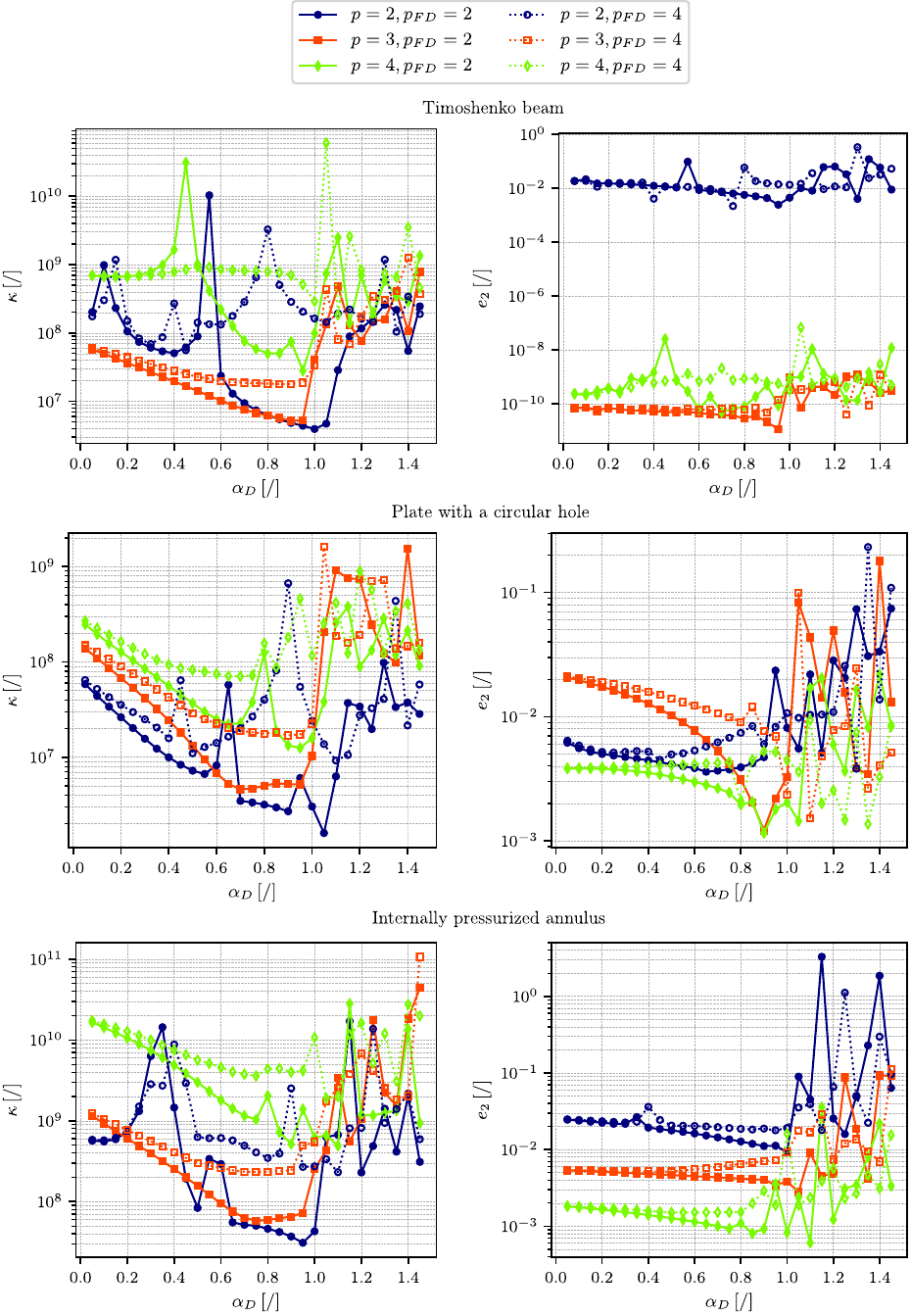} 
		\caption{Condition number $\kappa$ (top row) and a relative error $e_2$ (bottom row) as a function of $\alpha_D$ for Timoshenko beam (left), plate with a circular hole (center) and internally pressurized annulus (right) at different values of augmentation $p$ and FD order $p_{FD}$.}
		\label{fig: K_e2_ALL}
	\end{centering}
\end{figure}

In Figure \ref{fig: K_e2_ALL}, $\kappa$ and $e_2$ are plotted as a function of $\alpha_D$ for all cases introduced. The results were obtained with the same average node spacing of $h = 0.033$ m. Calculations were performed for values of $\alpha_D$ from $0.05$ to $1.95$. Plots show solutions for $\alpha_D$ only up to $1.45$ since the behavior at larger values is dominated by noise, and no coherent conclusions can be made. SNs that lie near the boundary may potentially fall outside the computational domain. To prevent this from happening and to keep FD stencils inside the domain, each SN is tested to see if it lies inside the domain. If not, $\alpha_D$ is set to $\alpha_{D, max} = 1$ for $p_{FD} = 2$ and $\alpha_{D, max} = 0.5$ for $p_{FD} = 4$. When comparing the values of $\alpha_D$ in Figure \ref{fig: K_e2_ALL}, the reader should keep in mind that at $\alpha_D = 1$, the stencil size for $p_{FD} = 2$ is $2h$, and for $p_{FD} = 4$, is $4h$, as defined in Figure \ref{fig: new_staff_1}. Again, no stabilization or other conditions are imposed on the boundary.

For the Timoshenko beam (Figure \ref{fig: K_e2_ALL}, left), it can be seen that for $p_{FD} = 2$, values of $\kappa$ experience some jumps but generally decrease when $\alpha_D \to 1$. For $\alpha_D \gtrsim 1$, the behavior becomes more unstable where values of $\kappa$ start to oscillate and increase in general. For results of $p_{FD} = 4$, a smaller reduction in $\kappa$ is observed when $\alpha_D \to 1$ or even no reduction for the case of $p=2$. The maximal change in values of $\kappa$ is of the order $\sim 10^{3}$. Similarly, a reduction in error can be observed where $\alpha_D \to 1$. As in the convergence study of the \textit{composed} approach presented in the previous section, one can also see that with $p=2, 3$, the exact solution is obtained. The maximum change in error is less than $\sim 10^{1}$.       

Fewer jumps can be observed for the plate with a circular hole (Figure \ref{fig: K_e2_ALL} center) where $\alpha_D \lesssim 1$. It can be seen that with a higher value of $p$, a higher value of $\kappa$ is obtained in general. As in the Timoshenko beam, $\kappa$ does not decrease as much for $p_{FD} = 4$ as for $p_{FD} = 2$ when $\alpha_D \to 1$. The maximal change that occurs in $\kappa$ is of the order $\sim 10^{2}$. Smooth transitions in error can be observed where $\alpha_D \lesssim 0.8$. For $p_{FD} = 4$ and $p=2, 4$, as $\alpha_D \to 1$, $e_2$ slightly increases. Here oscillatory behavior in all cases starts a little earlier, around $\alpha_D \approx 0.9$. The maximal change in error is of the order $\sim 10^{1}$. 

In Figure \ref{fig: K_e2_ALL} (right), similar behavior of the internally pressurized annulus can be seen. Again, the maximal change that occurs in $\kappa$ is of the order $\sim 10^{2}$. In this case, a very small reduction in $e_2$ can be noticed if $p_{FD} = 2$. When $p_{FD} = 4$, a slight increase occurs, especially at $p=3, 4$.   

For the relative error, it is known \cite{fasshauer_meshfree_2007} that it consists of interpolation error and conditioning of the problem. Interpolation error is for FD method expected to be decreasing with decreasing stencil size. Also, in the RBF-FD, it is known that evaluating operators in regions between collocation nodes induces error. From the results where the error generally decreases with increasing $\alpha_D$, we can claim that the conditioning of the problem governs relative error. One can also see that the jumps in $e_2$ (for example, Timoshenko beam with $p=4$) appear at the same $\alpha_D$ as jumps in $\kappa$.

To sum up, with the \textit{hybrid} approach in use, the size of $\alpha_D$ should be chosen in the stable regime, which has been shown to be when $\alpha_D \lesssim 1$. When $0.8 < \alpha_D < 1$, the smallest value in $\kappa$ is expected to be reached, so the method is expected to be most stable with the choice inside this interval. Similarly, the method is expected to be most accurate within this interval since $e_2$ should reach its minimum value.

\subsection{Effect of $\alpha_S$ on stability, error and convergence}\label{sec: eff_A_s}
In this section, the study of BCs stabilization is presented. The effect of the newly introduced parameter $\alpha_S$ is studied in terms of stability, accuracy, and convergence of the method. Discretization is performed with the \textit{hybrid} approach where $\alpha_D = 0.5$ and $p_{FD} = 2$. PHS order is set to $m=3$ and the order of augmentation to $p=2$. Results of condition number $\kappa(\alpha_S)$ and relative error $e_2(\alpha_S)$ are obtained for values of $0 < \alpha_S < 1.05$, where for each benchmark test, SNA is kept the same. The average node spacing was set to $h = 0.033$ m. 


\begin{figure}[H]
	\subfloat[]{%
		\includegraphics[width=0.5\textwidth]{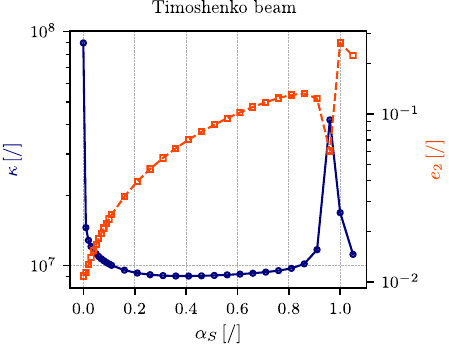}}
	\subfloat[\label{L2_err_H_SH_FD2} ]{%
		\includegraphics[width=0.5\textwidth]{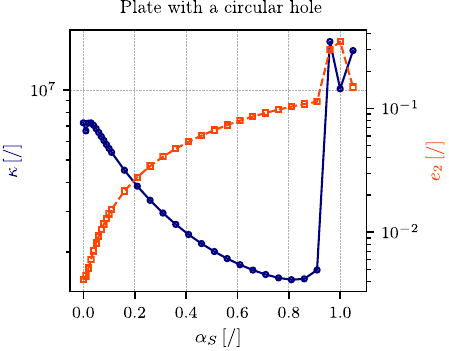}}\\
	\centering
	\subfloat[\label{L2_err_H_Sp0_FD2} ]{%
		\includegraphics[width=0.5\textwidth]{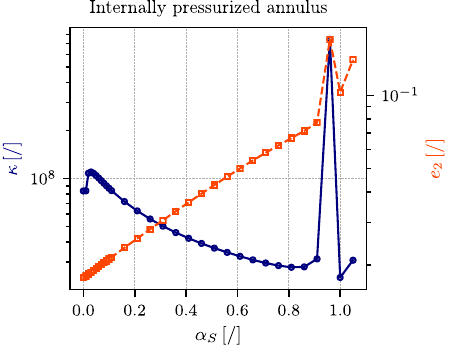}}\\
	\caption{\label{fig: skupna_a_S} Condition number $\kappa$, and a relative error $e_2$ as a function of $\alpha_S$ for (a) Timoshenko beam, (b) plate with a circular hole and (c) internally pressurized annulus, obtained with a \textit{hybrid} approach where $\alpha_D = 0.5$, $p_{FD} = 2$, $p=2$ and constant SNA with $h = 0.033$ m.}
\end{figure}

In Figure \ref{fig: skupna_a_S}, one can clearly see that with introducing boundary node shift ($\alpha_S \neq 0$), relative error starts to increase. In the Timoshenko beam, the relative error smoothly increases up to $\alpha_S = 0.9$. Then, when $\alpha_S \approx 1$, the behavior becomes erratic. A similar trend can be observed for the other two cases. On the other hand, the condition number $\kappa$ decreases with increasing $\alpha_S$. In the Timoshenko beam, $\kappa$ drops immediately for an already very small value of $\alpha_S$ and reaches its minimum at $\alpha_S \approx 0.4$. For the other two cases, $\kappa$ experiences some small jumps for a very small value of $\alpha_S$ and then decreases slowly until it reaches its minimum at $\alpha_S \approx 0.8$. This decreasing trend tells that stability is generally improved with $\alpha_S$ being introduced. However, the error increases since BCs are not met on the boundary. One can see that within the range of performed studies, $\kappa$ and $e_2$ do not change for more than one decade.
To stabilize the problem, $\alpha_S$ should be chosen as large as possible but not larger than $\alpha_S = 0.8$. If the need for stabilization is not so high, smaller values can be used for better accuracy. 

\begin{figure}
	\begin{centering}
		\includegraphics[width=1\textwidth]{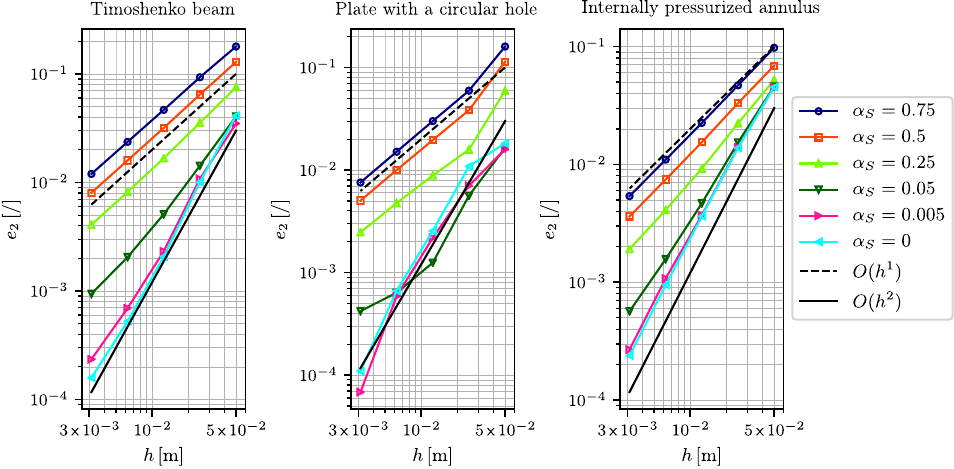}
		\caption{Convergence of \textit{hybrid} approach for different values of $\alpha_S$ for Timoshenko beam (left), plate with a circular hole (center) and internally pressurized annulus (right) where $\alpha_D = 0.5$, $p_{FD} = 2$ and $p=2$.}
		\label{fig: L2_err_New}
	\end{centering}
\end{figure}

Figure \ref{fig: L2_err_New} represents the convergence test for different values of $\alpha_S$. One can see similar behavior in all three cases. It can be seen that the convergence is already spoiled for very small values of $\alpha_S$. For $\alpha_S \gtrsim 0.25$, approximately the first order of convergence is obtained. 
If the need for stabilization is high, it is better to choose larger $\alpha_S$ since stability is increased and the convergence will not change significantly from $\alpha_S \approx 0.25$ up.

To sum up, the value of $\alpha_S$ should be chosen within the $0.25 \lesssim \alpha_S \lesssim 0.8$ range. Small values will not have much impact on stability and will quickly increase the error. Big values $\alpha_S \gtrsim 0.8$ will spoil the stability. 

\subsection{Effect of omission of corner nodes on the solution}
\begin{figure}[H]
	\begin{centering}
		\includegraphics[scale=1]{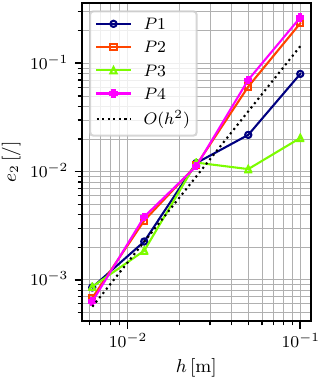}
		\caption{The convergence of extrapolated displacements at the corner nodes of a quarter of an annulus.}
		\label{fig:Corner}
	\end{centering}
\end{figure}

Due to reasons described in Section \ref{sec: Geom_dist}, the discretization nodes are not positioned in the corners of a domain. To study the effect of corner nodes omission, we test the convergence of the extrapolated solution in corners on a case of internally pressurized annulus. In Figure \ref{fig:Corner}, the error defined by equation \ref{eq: e2} (with N=1) is shown for corner nodes $P1 = [Ri,0], P2 = [0, Ri], P3 = [Ro, 0], P4 = [0, Ro]$ (see Figure \ref{fig: EL_ALL}), where the numerical solution is obtained by extrapolation. Results are obtained with the \textit{hybrid} approach, with $\alpha_D = 0.5, \alpha_S=0$, and the $2^{\text{nd}}$ order of augmentation.
From Figure \ref{fig:Corner}, we can see that the solution is converging with the same order as the proposed method, from where it can be seen that including corner nodes is not essential. If a critical area of interest lies in the corner, then additional node refinement can be applied in this region.

\section{Numerical performance in the elasto-plastic range} \label{sec: NUM_EX_PL}
This section studies the performance of introduced approaches on the internally pressurized annulus. Internal pressure load is applied incrementally as $\{p_{min}, p_{max}, \Delta p\} = \{8, 10.5, 0.1\}$ Pa. The initial yield stress is set to $\sigma_{y0} = 20$ Pa, and the hardening modulus is set to $H = 0$ Pa. This elasto-perfectly plastic response ($H = 0$ Pa) is known to be computationally much more challenging in terms of convergence than using non-zero hardening \cite{belytschko_nonlinear_2014}. Other values of used parameters are the same as presented in Section \ref{sec: IPA}. To reduce the computational time, and since the solution is radially symmetric, the geometry of observation is an annulus section with an angle of $\pi/6$. An example of SNA of the observed geometry is shown in \ref{app: Meshes}.

The reference solution (RS) was obtained by the finite element method (FEM). The well-established Abaqus program package \cite{smith_abaqusstandard_2009} was used. Discretization was done with 321920 linear quadrilateral elements where the size of the element edge $h_{FEM}$ was approximately $1/2$ of the minimal node distance $h$. The same geometry, load stepping and material parameters were used in obtaining RS.

\subsection{Analysis of the approaches without boundary stabilization}\label{sec: NoStab}
\begin{figure}[H]
	\begin{centering}
		\includegraphics[scale=1]{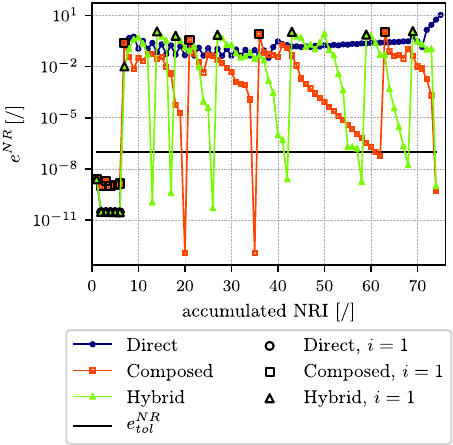}
		\caption{NRIA residual $e^{NR}$ for the first 75 accumulated NRIs for all proposed approaches obtained without BC stabilization.}
		\label{fig:NRres}
	\end{centering}
\end{figure}

With no BC stabilization imposed ($\alpha_S = 0$), none of the introduced approaches provide the correct solution. In Figure \ref{fig:NRres}, the residuals of NRIA ($e^{NR}$) are shown for the first 75 accumulated Newton-Raphson iterations (NRIs), where black contours represent the first residual ($i=1$) at some load increment $n$.

\begin{figure}[H]
	\subfloat[]{%
		\includegraphics[width=0.5\textwidth]{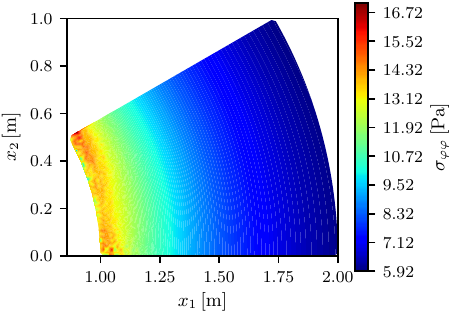}}
	\subfloat[ ]{%
		\includegraphics[width=0.5\textwidth]{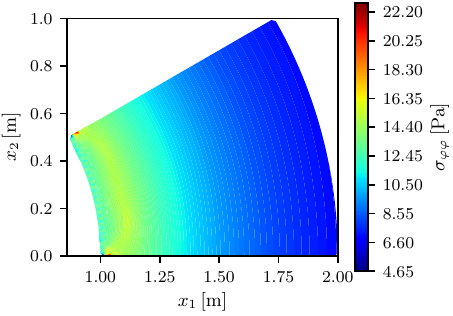}}
	\caption{\label{fig: NO_BCstab} Solution obtained without stabilization on the boundary; (a) last successfully converged solution $\sigma_{\varphi \varphi}(\bm{p})$ of the \textit{composed} approach, (b) $\sigma_{\varphi \varphi}(\bm{p})$ solution obtained with the \textit{hybrid} approach.}
\end{figure}

It can be seen that for the load increment $n \in \left\{1, 2, ..., 6\right\}$, the response is purely elastic. The convergence tolerance is set to $e_{tol}^{NR} = 10^{-7}$, and the maximum number of the NRIs is set to NRI$_{max} = 70$. Other parameters are set to $h = 0.0125$ m, $m=3$, $p=2$ and for the \textit{hybrid} approach $p_{FD} = 2$, $\alpha_{D} = 0.5$. One can see that the \textit{direct} approach does not converge in terms of NRIA. This happens since the stiffness matrix and the calculation of internal force are not consistent with each other. Similar inconsistency was applied in \cite{strnisa_meshless_2022}, where the direct iterative method was employed, and many global iterations were required to obtain the converged solution. Here, using NRIA, the proposed approach is found to be unsuitable for solving elasto-plastic problems. Using the \textit{composed} approach, NRIA converges successfully up to the $5^{\text{th}}$ plastic loading step. The solution of the hoop stress $\sigma_{\varphi \varphi}(\bm{p})$ is shown in Figure \ref{fig: NO_BCstab} (a). One can observe that in the plastic regime, the solution is oscillatory. With the \textit{hybrid} approach, NRIA converges to the last applied loading step. As shown in Figure \ref{fig: NO_BCstab} (b), the solution has no axisymmetric nature, as expected, but is generally smoother than the one obtained with the \textit{composed} approach. Since the \textit{direct} approach is found to be unusable for solving elasto-plastic problems, only the results obtained with the \textit{composed} and \textit{hybrid} approach are presented in the following.

\subsection{Composed approach performance}\label{sec: COM_app_elpl}
\begin{figure}[H]
	\begin{centering}
		\includegraphics[scale=1]{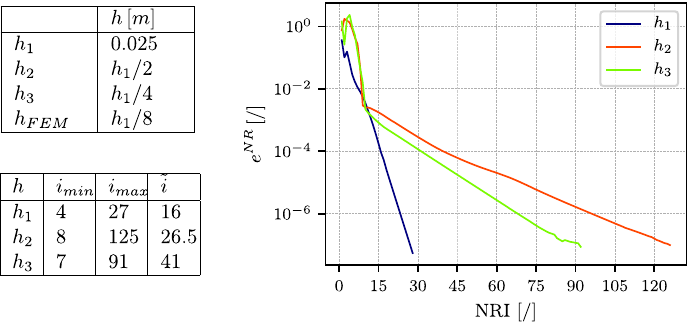}
		\caption{Average node spacings used in the study (top-left). Number of NRIs $i$ needed for reaching convergence tolerance $10^{-7}$ for three different SNA densities (bottom-left) and NRIA residual for the loading step with the slowest rate of convergence (right) obtained with the \textit{composed} approach.}
		\label{fig: NR_composed}
	\end{centering}
\end{figure}

In this study, all numerical parameters, except for $\alpha_S$, are the same as those in Section \ref{sec: NoStab}. Due to a high need for stabilization of BCs we set $\alpha_S= 0.5$. The proposed approach is studied on three different SNA densities, defined in Figure \ref{fig: NR_composed} (top-left).

Figure \ref{fig: NR_composed} (bottom-left) shows the number of NRIs needed for the convergence. The minimum number for the convergence in some loading step is denoted as $i_{min}$, the maximum number as $i_{max}$ and the median as $\tilde{i}$. Only the elasto-plastic loading steps are taken into account. One can see that convergence in terms of NRIA is relatively slow. It is observed that the number of NRIs is, on average, increasing when the plastic zone increases. It can be seen from the median values that, with increasing SNA density, the number of NRIs also increases.  

The loading step with the slowest convergence of NRIA is shown in Figure \ref{fig: NR_composed} (right). As in Figure \ref{fig:NRres}, the oscillatory behavior is still present that eventually smooths out with the increasing number of NRIs. 

\begin{figure}[H]
	\begin{centering}
		\includegraphics[width=1\textwidth]{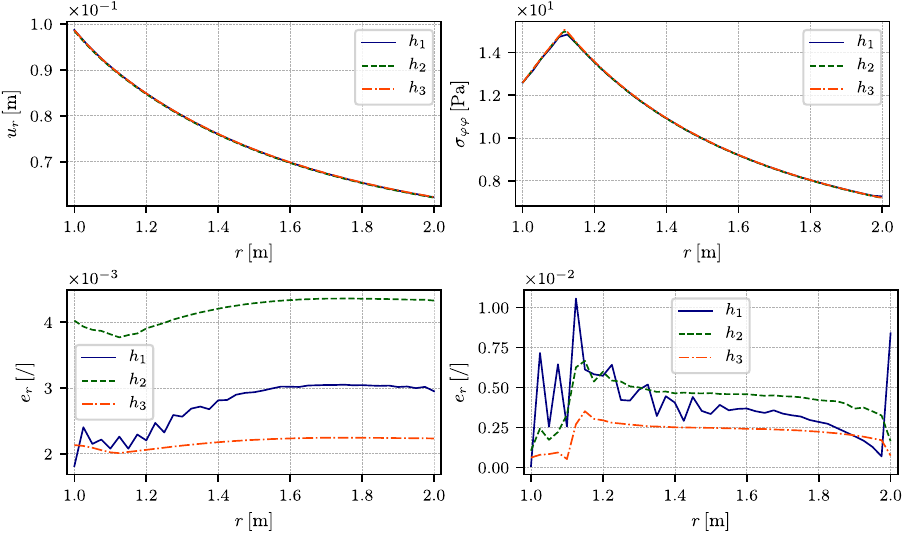}
		\caption{Results of $u_r$ and $\sigma_{\varphi \varphi}$ over the line for three different SNA densities and the RS (top) obtained with the \textit{composed} approach. Corresponding relative error (bottom).}
		\label{fig: composed_ALL}
	\end{centering}
\end{figure}

Solutions of radial displacement $u_r$ and hoop stress $\sigma_{\varphi \varphi}$ from the last loading step are presented in Figure \ref{fig: composed_ALL}. Results are shown over the line that splits the observed geometry in half in the radial direction. A good match is observed between an RS and the \textit{composed} approach solutions. Relative error over the line, shown on the bottom two plots of Figure \ref{fig: composed_ALL}, is defined as 
\begin{equation}
	e_r(r) = \frac{| y(r) - \hat{y}^{RS}(r) |}{|\hat{y}^{RS}(r)|},
\end{equation}
where $\hat{y}^{RS}$ represents the scalar value of the reference solution and $y$ of the \textit{composed} or \textit{hybrid} approach.

\begin{figure}[H]
	\begin{centering}
		\includegraphics[scale=1]{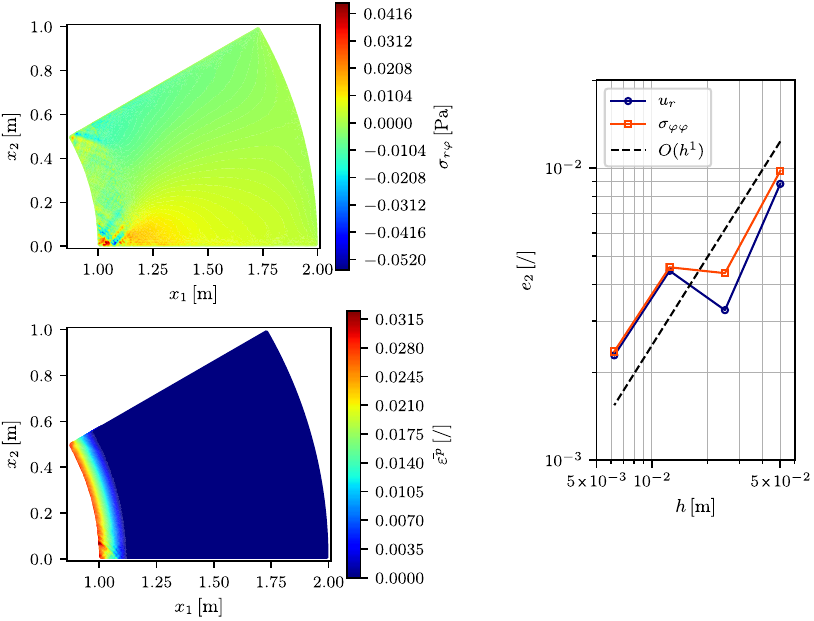}
		\caption{Results obtained with the \textit{composed} approach. Shear component of stress tensor $\sigma_{r \varphi}(\bm{p})$ (top-left), and accumulated plastic strain $\bar{\varepsilon}^p(\bm{p})$ (bottom-left). Convergence in terms of SNA density (right).}
		\label{fig: convergenca_composed_w_FV}
	\end{centering}
\end{figure}

Except for the case with $h_1$, the error in $u_r$ is smallest at the elasto-plastic front (EPF), which separates the region with $\bar{\varepsilon}^p = 0$ and $\bar{\varepsilon}^p \neq 0$. In the elastic region, it smoothly reaches a constant value. The error in $\sigma_{\varphi \varphi}$, on the other hand, experiences a jump at EPF and also reduces at the boundary. On the coarse SNA ($h_1$), it has an oscillatory shape. Oscillations can also be observed in the solution on a very dense SNA ($h=h_3$) in Figure \ref{fig: convergenca_composed_w_FV} (left). The shear component of a stress tensor $\sigma_{r \varphi}(\bm{p})$, which should be equal to zero, is oscillatory, especially in the plastic region. A similar effect can be observed in the accumulated plastic strain $\bar{\varepsilon}^p(\bm{p})$.
As shown in Section \ref{sec: eff_A_s}, using $p=2$ and $\alpha_S \gtrsim 0.25$, approximately one order in convergence is expected to be lost. Here, using $\alpha_S= 0.5$, similar observations can be made from the plot of convergence of $e_2$ in terms of $u_r$ and $\sigma_{\varphi \varphi}$ given in Figure \ref{fig: convergenca_composed_w_FV} (right). It is interesting to see that the error of $\sigma_{\varphi \varphi}$, which is associated with the derivative of $u_r$, converges at the same rate as that of $u_r$.

\subsection{Hybrid approach performance}
In this section performance of the \textit{hybrid} approach is presented. The parameters of the method, except for $\alpha_D$, are the same as in Section \ref{sec: COM_app_elpl}. In Section \ref{sec: A_D}, it is shown that the value of $\alpha_D$ should be chosen as $0.8 < \alpha_D < 1$ to stay in a stable regime and to obtain the best accuracy and stability. Within this range, $\kappa$ reaches a minimum, so NRIA is expected to converge at the highest rate. On the other hand, as shown in Section \ref{sec: Con_HY}, with smaller values of $\alpha_D$, the stability in terms of convergence is increased. The performance is thus studied for values of $\alpha_D = \left\{0.1, 0.3, 0.5, 0.8\right\}$. 

\begin{figure}
	\begin{centering}
		\includegraphics[scale=1]{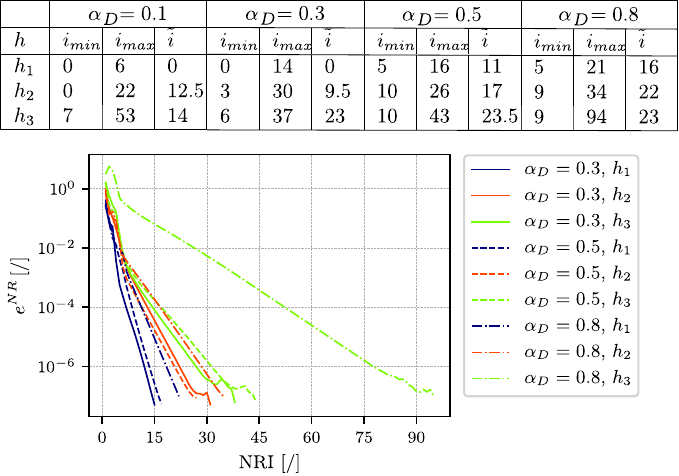}
		\caption{Number of NRIs $i$ needed for reaching convergence tolerance $10^{-7}$ using the \textit{hybrid} approach (top). NRIA residual for the loading step with the slowest convergence rate (bottom).}
		\label{fig: NR_hybrid}
	\end{centering}
\end{figure}

Figure \ref{fig: NR_hybrid} (top) shows the number of NRIs needed for the convergence at different SNA densities and different values of parameter $\alpha_D$. In cases where $i_{min} = 0$, the solution procedure fails. As shown, this happens with small values of $\alpha_D$ and coarse SNAs. With $\alpha_D = 0.5$ and $\alpha_D = 0.8$, the correct solution is also obtained on SNA with $h_1$. The number of NRIs $i$ is generally lower at the beginning, where a small number of SNs fall into the zone with $\bar{\varepsilon}^p \neq 0$ and then increase as this zone gets larger. It should also be emphasized that with decreasing $h$ and $\alpha_D$ the SNs, where the stresses are evaluated, are moving further from the boundary. This means that the EPF will reach SNs at a higher pressure load. Ignoring $\alpha_D = 0.1$, one can see that the difference in the median is getting smaller with increasing SNA density. It has also been observed that the oscillatory behavior of the residual of NRIA is generally less prominent at higher values of $\alpha_D$. 

In Figure \ref{fig: NR_hybrid} (bottom), the loading step with the slowest convergence of NRIA is shown. Comparing results with the $\textit{composed}$ approach in Figure \ref{fig: NR_composed} (right), one can see that fewer NRIs is needed to reach convergence tolerance. Since results obtained with $\alpha_D = 0.1$ are found to be incorrect, the solutions for this case are not shown. It is evident that with increasing $\alpha_D$ and SNA density, the number of NRIs is increasing. With $\alpha_D = 0.8$ and $h_3$, the $i_{max}$ is almost twice as much as with $\alpha_D = 0.5$. The elastic benchmark study suggested that the fastest convergence of NRIA would be achieved by selecting $\alpha_D$ within the region of the minimum condition number. This is not observed.

The first part of the paper shows that $\alpha_D$ should be smaller than $1$ or optimally $0.8 < \alpha_D < 1$. Here it is shown that $\alpha_D$ should not be too small $\alpha_D \gtrsim 0.3$ to have a stable convergence in terms of NRIA and not too big $\alpha_D \lesssim 0.8$ to obtain as few NRIs as possible. So, in the following, $\alpha_D = 0.5$ was chosen to obtain results.

\begin{figure}[H]
	\begin{centering}
		\includegraphics[width=1\textwidth]{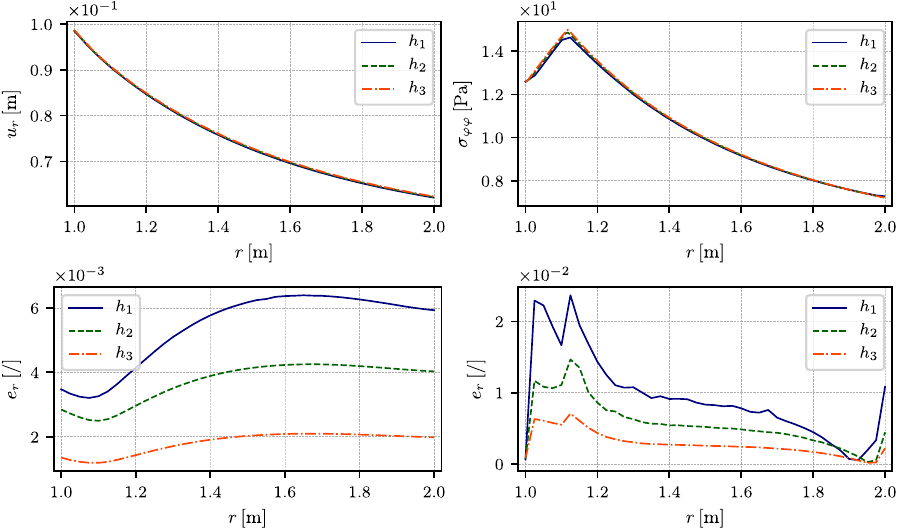}
		\caption{Results of $u_r$ and $\sigma_{\varphi \varphi}$ over the line for three different SNA densities with $\alpha_D= 0.5$ and the RS (top) obtained with the \textit{hybrid} approach. Corresponding relative error (bottom).}
		\label{fig: hybrid_ALL2}
	\end{centering}
\end{figure}

\begin{figure}[H]
	\begin{centering}
		\includegraphics[scale=1]{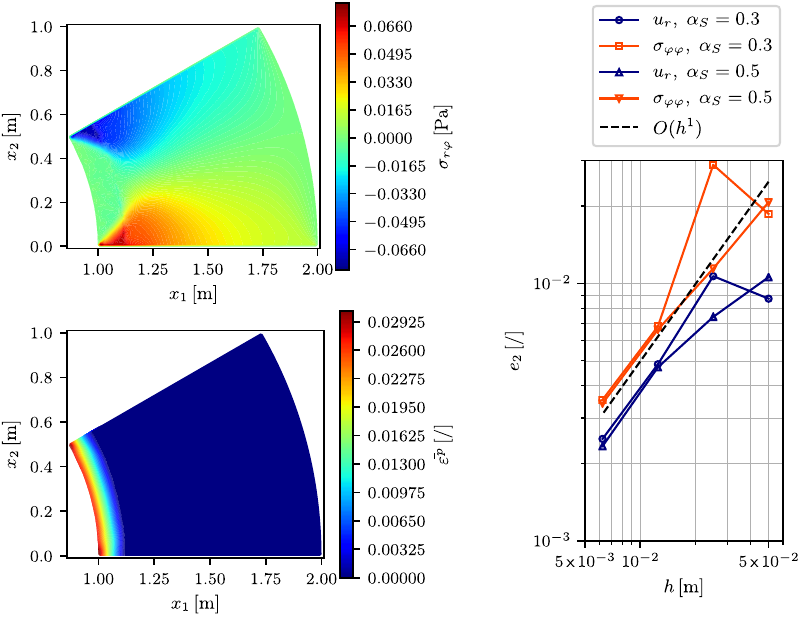}
		\caption{Results obtained with the \textit{hybrid} approach. Shear component of stress tensor $\sigma_{r \varphi}(\bm{p})$ (top-left), and accumulated plastic strain $\bar{\varepsilon}^p(\bm{p})$ (bottom-left). Convergence in terms of SNA density (right).}
		\label{fig: convergenca_hybrid_w_FV}
	\end{centering}
\end{figure}

As in the \textit{composed} approach, the results of $u_r$ and $\sigma_{\varphi \varphi}$ are compared with an RS and shown in Figure \ref{fig: hybrid_ALL2}. A good match can be observed. Compared with the \textit{composed} approach, a relative error is much smoother, but in size, it is approximately two times larger in terms of $\sigma_{\varphi \varphi}$ and approximately the same in terms of $u_r$. A larger error in stress can also be observed in the shear stress component $\sigma_{r\varphi}(\bm{p})$ presented in Figure \ref{fig: convergenca_hybrid_w_FV} (top-left). The accumulated plastic strain is shown in Figure \ref{fig: convergenca_hybrid_w_FV} (bottom-left). From both field values, it
can be seen that the oscillatory behavior is not present as in the \textit{composed} approach. 


Similarly, as in the \textit{composed} approach and elastic cases, when $\alpha_S = 0.5$, the first-order convergence is obtained, as shown in Figure \ref{fig: convergenca_hybrid_w_FV} (right) for $\alpha_D=0.3$ and $\alpha_D=0.5$. It can be seen that in terms of accuracy, the \textit{composed} approach performs better.

\subsection{Comparison of composed and hybrid approach in terms of approximation-induced discontinuities}
This section presents a study on the approximation-induced discontinuities in a solution field for a \textit{composed} and \textit{hybrid} approach. Suppose Voronoi tessellation is performed in such a way that each CN represents a seed of the corresponding Voronoi cell. It is well known \cite{fasshauer_meshfree_2007} that when using augmented RBF-FD, the solution on the boundary of two contact Voronoi cells depends on which neighboring support domain is chosen to interpolate the solution on the boundary. In other words, approximation-induced discontinuity (AID) in interpolated solution is obtained on the border of contact Voronoi cells. It is also known that AID is decreasing as $h \to 0$ \cite{tominec_least_2021}.

Lets now observe the inner force vector as $\ve{f}^{int} = \nabla \cdot \bm\sigma = \nabla \cdot \left(\textbf{\textsf{D}} \bm\varepsilon \right)=\nabla \cdot \left(\textbf{\textsf{D}} \nabla^s \ve{u} \right)$. In the stress field, two contributions of discontinuity are present. One arises from the AID and the other from the material property in $\textbf{\textsf{D}}$. The divergence operator is with the \textit{composed} approach applied in an RBF-FD manner where it is generally assumed that the field on which the operator acts is continuous, which is not the case here. With the \textit{hybrid} approach, the divergence operator is evaluated inside (assuming $\alpha_D$ is small enough) the Voronoi cells, so the derivation over AIDs is overcome here. This affects the displacement field solution and its interpolation on the boundary of the Voronoi cells. 

\begin{figure}[H]
	\begin{centering}
		\includegraphics[scale=1]{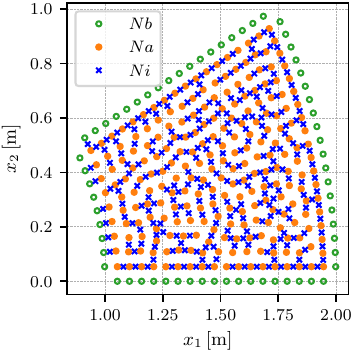}
		\caption{Splitting nodes $Ni$ positioning, where $Na$ and $Nb$ represent inner and boundary nodes, respectively.}
		\label{fig: mreza2}
	\end{centering}
\end{figure}

\begin{figure}[H]
	\begin{centering}
		\includegraphics[scale=0.62]{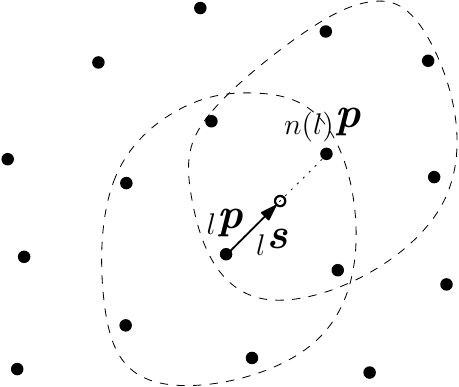}
		\caption{Two overlapping domains where the center node ${}_{l}\ve{p}$ of domain ${}_{l}\Omega$ has its nearest neighbor ${}_{n(l)}\ve{p}$. With vector ${}_{l}\ve{s}$, the splitting point (hollow circle) is positioned halfway between points ${}_{l}\ve{p}$ and ${}_{n(l)}\ve{p}$.}
		\label{fig: Nezveznost}
	\end{centering}
\end{figure}

The following study on AIDs is performed for the \textit{composed} and the \textit{hybrid} approach. Instead of performing Voronoi tessellation and determining Voronoi cell boundaries, a set of splitting nodes is defined, as shown in Figure \ref{fig: mreza2}. Each splitting node is positioned between the two nearest ${}_{l}\Omega$-s, as shown in Figure \ref{fig: Nezveznost}, with a hollow circle. Vector ${}_{l}\ve{s}$ is given as
\begin{equation}
	{}_{l}\ve{s} = ({}_{n(l)}\ve{p} - {}_{l}\ve{p})/2,
\end{equation}
where function $n(l)$ defines the nearest neighbor to ${}_{l}\ve{p}$. If during computation of $n(l)$, where $l$ runs over inner nodes $Na$, the algorithm runs into position $n(l+k) = l$, where $n(l) = l+k $ was already determined for some $k$, next (second, third, or so on) nearest neighbor is prescribed. This procedure ensures that each splitting node has a different position and that the number of them is $Ni = Na$. The value of AID is assessed here with the relative error defined as 
\begin{equation}
	e_{int} = \sqrt{\frac{\sum_{l=1}^{Na} \left\| {}_{l}\ve{u} \left( {}_{l}\ve{p} + {}_{l}\ve{s} \right)  -  {}_{n(l)}\ve{u} \left( {}_{n(l)}\ve{p} - {}_{l}\ve{s} \right) \right\|^2 }{\sum_{l=1}^{Na} \left\| {}_{l}\ve{u} \left( {}_{l}\ve{p} + {}_{l}\ve{s} \right) \right\|^2}},
\end{equation}
where the sum runs only over inner nodes $Na$.

\begin{figure}[H]
	\begin{centering}
		\includegraphics[scale=1]{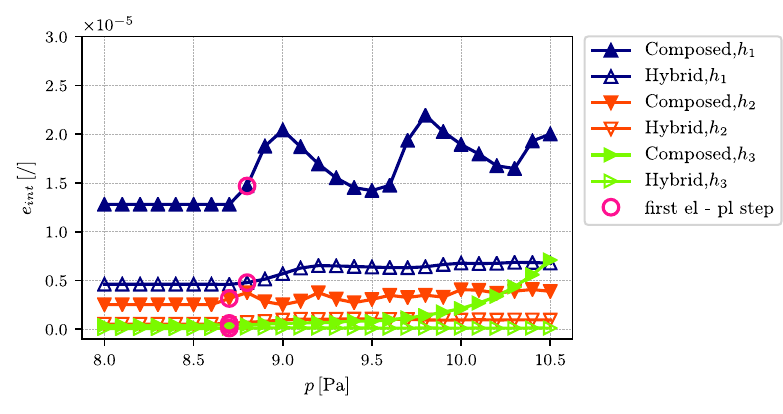}
		\caption{Relative error $e_{int}$ as a function of incremental pressure load for \textit{composed} and \textit{hybrid} approach at different SNA densities. Circles emphasize the first elasto-plastic loading step.}
		\label{fig: error_int2}
	\end{centering}
\end{figure}

Figure \ref{fig: error_int2} shows $e_{int}$ for \textit{composed} and \textit{hybrid} approaches as a function of incremental pressure load. Solution field of displacement $\ve{u}$ has been obtained using $m=3, p=2, p_{FD} = 2, \alpha_S= 0.5$ and $\alpha_D = 0.5$. The first elasto-plastic steps are denoted with circles. Within elastic load stepping, $e_{int}$ stays constant for all cases. It is evident that $e_{int}$ is for particular $h$ larger for the \textit{composed} approach. $e_{int}$ starts to increase with increasing load once the response is elasto-plastic. Using the \textit{hybrid} approach $e_{int}$ increases only at the beginning of the elasto-plastic response. It remains constant, unlike the \textit{composed} approach where $e_{int}$ is growing further.

\subsection{Comparison of hybrid approach with FEM}
For additional verification of the proposed \textit{hybrid} approach, the case of a plate with a hole introduced in Section \ref{sec: PWH} is studied in an elasto-plastic regime within a plane strain approximation and compared with FEM. In addition to elastic material parameters defined in Figure \ref{fig: EL_ALL}, the initial yield stress is set to $\sigma_{y0} = 0.1$ Pa. Hardening modulus is $H = E/4 = 0.25$ Pa. Traction load, derived from equation \ref{eq: PLW_sol}, is applied incrementally where load amplitude $\sigma_{\infty}$ is applied as $\left\{\sigma_{\infty, min}, \sigma_{\infty, max}, \Delta \sigma_{\infty}\right\}= \left\{0,0.1,0.01\right\}$ Pa. The parameters of the hybrid method are set to $\alpha_D = \alpha_S = 0.5$, and convergence tolerance is set to $e_{tol}^{NR} = 10^{-7}$. FEM solution is obtained with 6-node quadratic finite elements as in the elastic regime. All other parameters are the same as in the \textit{hybrid} approach.

\begin{figure}[H]
	\begin{centering}
		\includegraphics[scale=1]{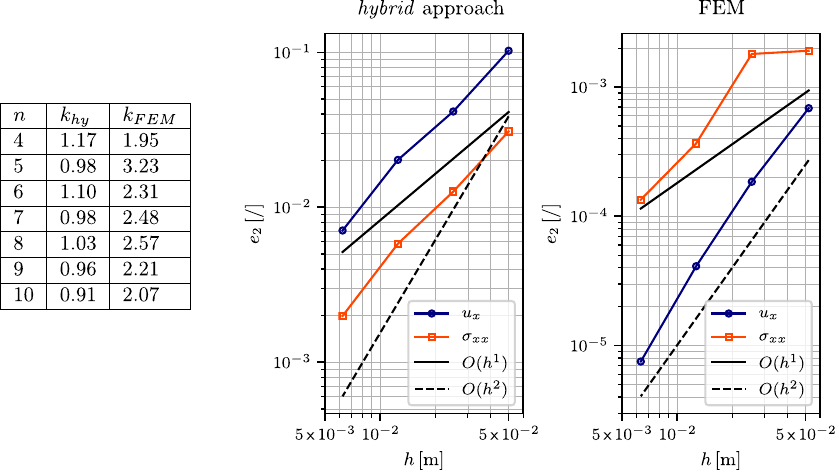}
		\caption{Measured averaged parameter $k$ for \textit{hybrid} approach $k_{hy}$ and FEM $k_{FEM}$ over plastic load increments $n$ (left). Convergence of the \textit{hybrid} approach (center) and convergence of FEM (right).}
		\label{fig: PWH_NR_CON}
	\end{centering}
\end{figure}

Using the NRIA for solving the system of equations, the error (in our case residual) is expected to converge by the following relation:
\begin{equation}
	\rho_{i+1} = A_0 \rho_{i}^k,	
\end{equation}
where $i$ is the iteration index, $A_0$ a constant and $k$ is the parameter that defines the order of convergence. Under certain conditions, it can be shown that NRIA converges with the $2^{\text{nd}}$ order ($k=2$) \cite{ortega_iterative_1970}. To test the convergence of the \textit{hybrid} approach and FEM, parameter $k$ is studied in terms of the largest residual force defined as $\rho= \text{max}\left(\left\{r_{x,i}, r_{y,i}\right\}; i \in \left\{1,N\right\}\right)$. Averaged results within each plastic load increment $n \in \left(4,10\right)$ of parameter $k$, measured for \textit{hybrid} approach ($k_{hy}$) and FEM ($k_{FEM}$), are shown in Figure \ref{fig: PWH_NR_CON} (left).
One can see that parameter $k_{FEM} \sim 2$. $k_{hy} \sim 1$ indicates that the \textit{hybrid} approach does not fulfill some of the conditions required to achieve $2^{\text{nd}}$ order convergence.

As in the previous case of internally pressurized annulus, the convergence in terms of node refinement is performed next. In Figure \ref{fig: PWH_NR_CON} (center) and (right), self-convergences, in terms of $u_x$ and $\sigma_{xx}$, are shown for the \text{hybrid} approach and FEM, respectively. The reference solution is obtained by employing half the node spacing used in the finest solution plotted. The convergence order of the \text{hybrid} approach is approximately $1.5$. Again, interestingly, displacement and stress values converge at the same rate, and additionally the error is smaller in stress. FEM results of displacements converge with the $2^{\text{nd}}$ order, and stresses converge with the $1^{\text{st}}$ order. Compared with the \textit{hybrid} approach, the error in FEM is smaller by about $10^2$, and stresses experience lower order of convergence.   

\begin{figure}[H]
	\subfloat[\textit{hybrid} approach \label{Svm_160_OURS}]{%
		\includegraphics[width=0.5\textwidth]{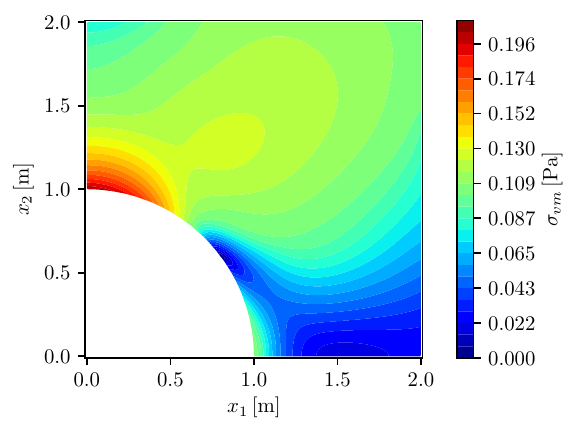}}
	\subfloat[FEM \label{Svm_160_FEM}]{%
		\includegraphics[width=0.5\textwidth]{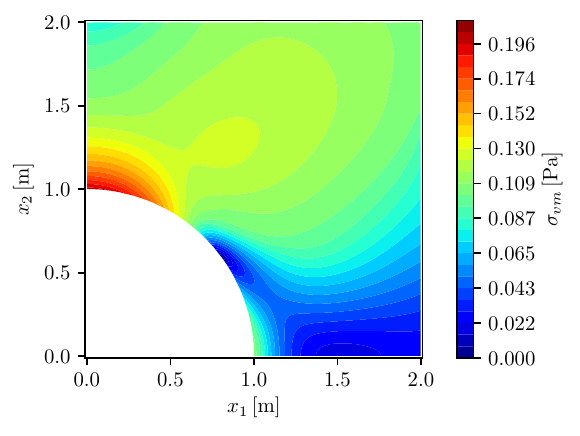}}\\
	\caption{\label{fig:Svm_160} Von Misses stress solution obtained with \textit{hybrid} approach (a), and FEM (b).}
\end{figure}

\begin{figure}[H]
	\subfloat[\textit{hybrid} approach \label{PEEQ_160_OURS}]{%
		\includegraphics[width=0.5\textwidth]{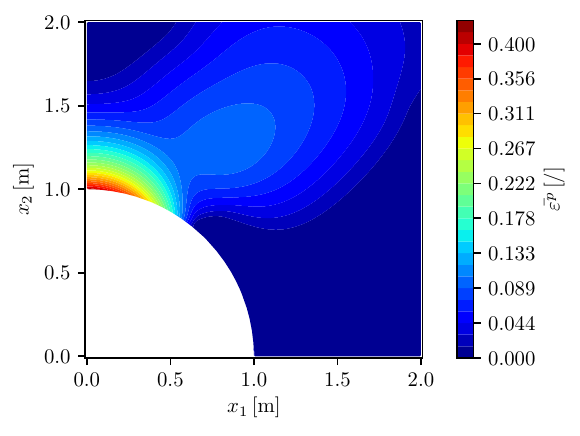}}
	\subfloat[FEM \label{PEEQ_160_FEM}]{%
		\includegraphics[width=0.5\textwidth]{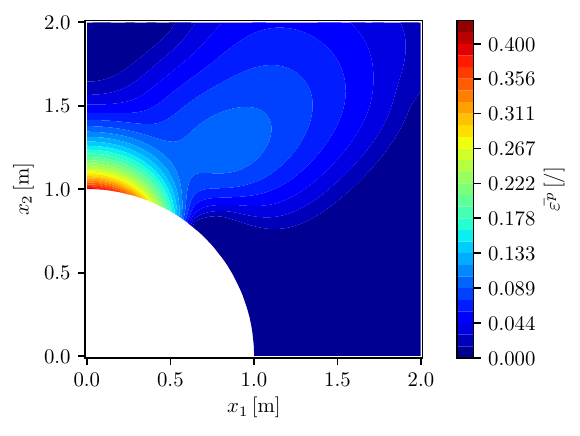}}\\
	\caption{\label{fig:PEEQ_160} Accumulated plastic strain solution obtained with \textit{hybrid} approach (a), and FEM (b).}
\end{figure}

Figures \ref{fig:Svm_160} and \ref{fig:PEEQ_160} show von Mises stress and accumulated plastic strain over the field at the final loading step. Solutions are plotted for a case obtained on the finest discretization used as a reference solution in the convergence study. Small changes with FEM solutions can be seen. The reason can be a different position of enforcing BCs since $\alpha_S > 0$ in the \textit{hybrid} approach.

\section{Conclusions}\label{sec: Conc}
This work introduced three approaches (\textit{direct}, \textit{composed} and \textit{hybrid}) for solving the elasto-plastic problems using the strong form meshless method.

The proposed approaches are studied on elastic benchmarks initially. The convergence is first checked on the \textit{composed} approach, where it is found that the order of augmentation governs the convergence order as it was previously confirmed for the \textit{direct} approach \cite{slak_adaptive_2020}. For the \textit{hybrid} approach, it has been shown that within $0.8 < \alpha_D < 1$, the method appears to be the most stable and accurate since the condition number and error reach the lowest value in that range. Due to similar trends of condition number and relative error, it is observed that the condition number governs the error. It was found that the convergence order is governed by augmentation, even when the order of FD is lower than the order of augmentation. It was also shown that with small values of $\alpha_D$, the convergence is uniform; with large values, the convergence is oscillatory and spoiled.

It was found that the method's stability is increased by introducing the parameter $\alpha_S$. Since the BCs are not evaluated exactly on the boundary, the accuracy deteriorates with an increase in $\alpha_S$. To stabilize BCs as much as possible, $\alpha_S$ should be chosen close to the best value of $\alpha_S \approx 0.8$. It was found that the order of convergence decreases as $\alpha_S$ increases. If $\alpha_S \gtrsim 0.25$, approximately the first order of convergence is obtained.


In the elasto-plastic case, performed after elastic ones, BC stabilization was shown to be essential. Without it, none of the introduced approaches provides a decent solution. Due to the inconsistent evaluation of the stiffness matrix and internal force in the \textit{direct} approach, it does not converge in NRIA. It is found to be unsuitable for solving elasto-plastic problems.
The \textit{composed} approach was found to be capable of solving elasto-plastic problems even more accurately than the \textit{hybrid} approach but suffers from oscillatory solutions. These oscillations are accumulated during loading causing slow convergence of NRIA. On average more NRIs are needed to converge compared to the \textit{hybrid} approach. With the \textit{hybrid} approach, the problem with oscillatory solution fields is overcome. Compared to the other two approaches, the downside of the \textit{hybrid} is four times (for $2^{\text{nd}}$ order FD) more nodes where stresses are evaluated. To obtain a stable convergence in terms of NRIA and to reach a minimum number of NRIs, $\alpha_D$ should be chosen as $0.3 \lesssim \alpha_D \lesssim 0.8$. 
In terms of approximation-induced discontinuity (AID), it was shown that the value of AID does not change much during incremental loading in the \textit{hybrid} approach. In contrast, with a \textit{composed} approach, the value of AID increases with load.

An additional comparison with FEM revealed that NRIA in FEM converges as expected with the $2^{\text{nd}}$ order. The \textit{hybrid} approach, on the other hand, does not converge with the same order in NRIA, indicating that some of the conditions for the $2^{\text{nd}}$ order convergence are not satisfied. Comparing convergence in node refinement showed that FEM is more accurate, but the benefit of the \textit{hybrid} approach is that stresses converge with the same order as displacements, which is not the case in FEM.

In future, the most promising \textit{hybrid} approach will be expanded to study more complex non-linear material models and will be applied to three-dimensional problems.

\section*{Acknowledgments}
Funding for this research is provided by the Slovenian Grant Agency (ARRS) within the framework of Young Researcher program, projects Z2-2640, J2-1317 and program group P2-0162. 


\section*{Statements and declarations}
\textbf{Conflict of interest}
The authors have no conflicts of interest to declare that are relevant to the content of this paper.

\appendix
\section{Hybrid approach discretization}\label{app: Operator_der}

\begin{figure}[!h]
	\begin{centering}
		\includegraphics[scale=1]{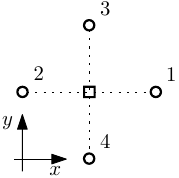}
		\caption{Enumeration of secondary nodes on a $2^{\text{nd}}$-order FD stencil.}
		\label{fig: alpha_D_enum}
	\end{centering}
\end{figure}

In the \textit{hybrid} approach, the coefficients of a gradient of a vector field have to be determined first. A similar derivation as in section \ref{sec: App_DO} gives
\begin{equation}
	\mathcal{L} {}_{l}\bm{y}(\ve{p})_{\xi \varsigma} \approx \displaystyle\sum_{\zeta=1}^{n_{d}} \displaystyle\sum_{j=1}^{{}_{l}N + M} {}_{l}\gamma_{j,\zeta} \displaystyle\sum_{\chi=1}^{n_d} \displaystyle\sum_{i=1}^{{}_{l}N+M} {}_{l}A_{ij,\chi \zeta}^{-1} \: \mathcal{L}_{\xi \varsigma \chi} \: {}_{l}\Psi_{i}(\ve{p}) = \displaystyle\sum_{\zeta=1}^{n_{d}} \displaystyle\sum_{j=1}^{{}_{l}N + M} {}_{l}\gamma_{j,\zeta} \: {}_{l}\mathscr{G}_{j, \zeta, \xi \varsigma} (\ve{p}),
\end{equation}
with 
\begin{equation}\label{eq: GOFVC}
	{}_{l}\mathscr{G}_{j, \zeta, \xi \varsigma} (\ve{p}) = \displaystyle\sum_{\chi=1}^{n_d} \displaystyle\sum_{i=1}^{{}_{l}N+M} {}_{l}A_{ij,\chi \zeta}^{-1} \: \mathcal{L}_{\xi \varsigma \chi} \: {}_{l}\Psi_{i}(\ve{p}).
\end{equation}
Next, coefficients for symmetric gradient operator $\nabla^s$ are defined. From equation (\ref{eq: GOFVC}) and definition of $\nabla^s$, it follows
\begin{equation}\label{eq: nabla_s}
	{}_{l,f}\mathscr{N}_{i, \zeta, \xi \varsigma} = \frac{{}_{l,f}\mathscr{G}_{i, \zeta, \xi \varsigma} + {}_{l,f}\mathscr{G}_{i, \zeta, \varsigma \xi}}{2}, \quad i \in [1, {}_{l}N+m], \quad \zeta,\xi, \varsigma \in [1, n_d], \quad l \in [1, N],
\end{equation}
where index $f \in \left[1,n_{dp}\right]$ runs over $n_{dp}$ SNs where the coefficients are evaluated. For the $2^{\text{nd}}$-order FD stencil in 2D, these are enumerated as shown in Figure \ref{fig: alpha_D_enum}. A similar procedure is done for the $4^{\text{th}}$-order FD stencil.
Next, multiplication with $\textbf{\textsf{D}}$ is performed. Since $\textbf{\textsf{D}}$ is defined in Voight notation and has in general dimension of $6 \times 6$ and for plane problems $4 \times 4$, mapping $m(\xi, \varsigma) = M_{\xi \varsigma}$ is introduced, which returns a value from the matrix
\begin{equation}\label{eq: map_M}
	M_{\xi \varsigma} = 
	\begin{bmatrix}
		1 & 4 & 5 \\
		4 & 2 & 6 \\
		5 & 6 & 3
	\end{bmatrix}.
\end{equation}
The coefficients ${}_{l,f}\mathscr{N}_{i, \zeta, \xi \varsigma}$ can thus be mapped to ${}_{l,f}\mathscr{N}_{i, \zeta, p}$. Multiplication is performed as
\begin{equation}
	{}_{l,f}\mathscr{M}_{i, \zeta, o} = \sum_{p=1}^{6} {}_{l,f}\textbf{\textsf{D}}_{i,o \, p} \:\: {}_{l,f}\mathscr{N}_{i, \zeta, p},
\end{equation}
to obtain coefficients of $(\textbf{\textsf{D}} \nabla^s)$. For clarity of presentation (without additional mappings), coefficients for $\nabla \cdot (\textbf{\textsf{D}} \nabla^s)$ for a 2D case in Cartesian coordinates are expressed as

\begin{equation}
	\begin{aligned}
		&{}_{l}\mathscr{D}_{i, \zeta, 1} =
		\frac{{}_{l,1}\mathscr{M}_{i, \zeta, m(1,1)} - {}_{l,2}\mathscr{M}_{i, \zeta, m(1,1)}}{2 \delta x} + \frac{{}_{l,3}\mathscr{M}_{i, \zeta, m(1,2)} - {}_{l,4}\mathscr{M}_{i, \zeta, m(1,2)}}{2 \delta y}\\
		&{}_{l}\mathscr{D}_{i, \zeta, 2} =
		\frac{{}_{l,1}\mathscr{M}_{i, \zeta, m(2,1)} - {}_{l,2}\mathscr{M}_{i, \zeta, m(2,1)}}{2 \delta x} + \frac{{}_{l,3}\mathscr{M}_{i, \zeta, m(2,2)} - {}_{l,4}\mathscr{M}_{i, \zeta, m(2,2)}}{2 \delta y}\\
	\end{aligned},
\end{equation}
where distances $\delta x = \delta y$ are defined in section \ref{sec: HybridApp} as $\delta x = \alpha_{D} \, h$.    
Equation (\ref{eq: MAIN_NR}) can be now for a single collocation node expressed as
\begin{equation}
	\sum_{\chi=1}^{n_d} \sum_{i=1}^{{}_{l}N} {}_{l}\mathscr{D}_{i, \zeta, \chi} \:\: {}_{l}\delta u_{i, \chi} = - {}_{l}r_{\zeta}.
\end{equation}
Internal force is expressed in the same manner as

\begin{equation}
	\begin{aligned}
		&{}_{l}f_{1} =  \frac{{}_{l,1}\sigma_{1} - {}_{l,2}\sigma_{1}}{2 \delta x} + \frac{{}_{l,3}\sigma_{4} - {}_{l,4}\sigma_{4}}{2 \delta y}\\
		&{}_{l}f_{2} =  \frac{{}_{l,1}\sigma_{4} - {}_{l,2}\sigma_{4}}{2 \delta x} + \frac{{}_{l,3}\sigma_{2} - {}_{l,4}\sigma_{2}}{2 \delta y}\\
	\end{aligned},
\end{equation}
where the stress values are obtained in the same manner as \ref{eq: map_M}.

\section{Examples of geometry discretization}\label{app: Meshes}
\begin{figure}[H]
	\begin{centering}
		\includegraphics[width=1\textwidth]{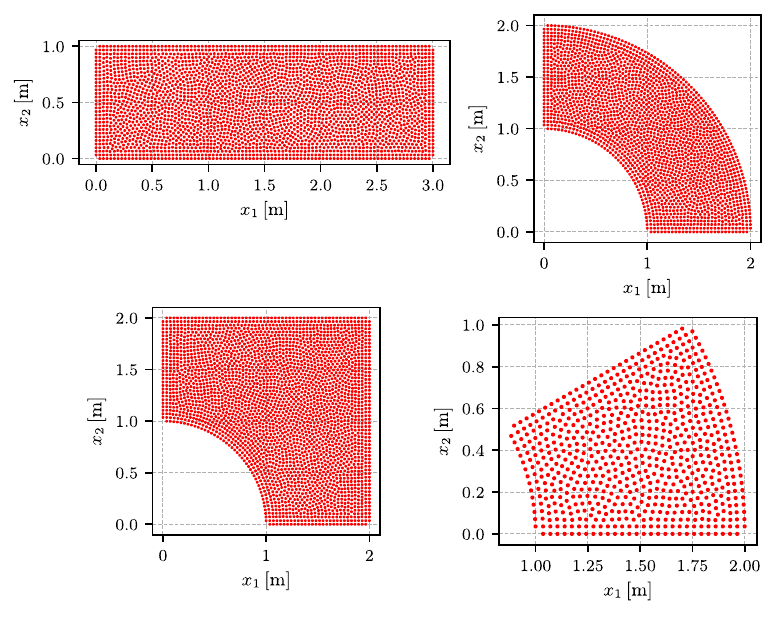}
		\caption{Examples of SNAs for Timoshenko beam (top-left), plate with a circular hole (bottom-left), internally pressurized annulus; elastic response (top-right), elasto-plastic response (bottom-right).}
		\label{fig: mrezenje_skupna}
	\end{centering}
\end{figure}

Examples of geometry discretization, where $h=0.033$ m, are presented in Figure \ref{fig: mrezenje_skupna}.


\bibliography{Article_1}

\end{document}